\documentclass[lettersize,journal]{IEEEtran}
\usepackage{amsmath,amsfonts}
\usepackage{algorithmic}
\usepackage{algorithm}
\usepackage{array}
\usepackage[caption=false,font=normalsize,labelfont=sf,textfont=sf]{subfig}
\usepackage{textcomp}
\usepackage{stfloats}
\usepackage{url}
\usepackage{verbatim}
\usepackage{graphicx}
\usepackage{cite}

\usepackage{mathtools}
\mathtoolsset{showonlyrefs}
\usepackage{amssymb}
\usepackage{amsthm}
\usepackage{xcolor}
\usepackage[shortlabels]{enumitem}
\usepackage{stmaryrd}
\usepackage{booktabs}

\newcommand{\re}[1]{{#1}}
\usepackage{cancel}

\normalfont

\renewcommand{\SS}{\ensuremath{\mathbb{S}}}

\newcommand{\R}{\ensuremath{\mathbb{R}}}

\newcommand{\T}{\mathcal{T}}
\newcommand{\XX}{\mathbb{X}}
\newcommand{\YY}{\mathbb{Y}}

\newcommand{\TT}{\ensuremath{\mathbb{T}}}

\newcommand{\supp}{\textnormal{supp}}

\newcommand{\p}{\mathcal{P}}

\newcommand{\dx}{\,\mathrm{d}}

\newcommand{\FS}{\mathfrak{S}}
\newcommand{\FX}{\mathfrak{X}}
\newcommand{\FY}{\mathfrak{Y}}

\newcommand{\Pio}{\Pi_\mathrm{o}}
\newcommand{\Gammas}{\Gamma_\sigma}
\newcommand{\geo}{\mathrm{g}}
\newcommand{\red}{\mathrm{r}}
\newcommand{\sym}{\mathrm{sym}}

\DeclareMathOperator*{\id}{id}

\DeclareMathOperator*{\argmin}{argmin}

\DeclareMathOperator{\GW}{GW}

\DeclareMathOperator{\LGW}{LGW}
\DeclareMathOperator{\gLGW}{gLGW}

\DeclareMathOperator{\LOT}{LOT}
\DeclareMathOperator{\gLOT}{gLOT}

\DeclareMathOperator{\Tan}{Tan}

\newtheorem{theorem}{Theorem}[section]
\newtheorem{lemma}[theorem]{Lemma}
\newtheorem{remark}[theorem]{Remark}

\newtheorem{example}[theorem]{Example}

\newtheorem{proposition}[theorem]{Proposition}

\allowdisplaybreaks

\begin{document}

\title{On a linear Gromov--Wasserstein distance}

\author{Florian Beier,
  Robert Beinert,
  Gabriele Steidl%
  \thanks{F. Beier is with the Institute of Mathematics,
	Technische Universit\"at Berlin, Stra\ss{}e des 17. Juni 136,
        10623 Berlin, Germany.}%
  \thanks{R. Beinert is with the Institute of Mathematics,
	Technische Universit\"at Berlin, Stra\ss{}e des 17. Juni 136,
        10623 Berlin, Germany.}%
  \thanks{G. Steidl is with the Institute of Mathematics,
	Technische Universit\"at Berlin, Stra\ss{}e des 17. Juni 136,
        10623 Berlin, Germany.}%
}

\markboth{Journal of \LaTeX\ Class Files,~Vol.~XX, No.~X, XXX}%
{Shell \MakeLowercase{\textit{et al.}}: A Sample Article Using IEEEtran.cls for IEEE Journals}

\IEEEpubid{0000--0000/00\$00.00~\copyright~2021 IEEE}

\maketitle

\begin{abstract}
Gromov--Wasserstein distances 
are generalization of 
Wasserstein distances,
which are invariant under 
distance preserving transformations.
Although a simplified version of optimal transport 
in Wasserstein spaces, called linear optimal transport (LOT),
was successfully used in \re{practice},
there does not exist a \re{notion}
of linear Gromov--Wasserstein distances so far.
In this paper, we propose a definition of linear Gromov--Wasserstein distances.
We motivate our approach by a generalized LOT model, 
which is based on barycentric projection maps of transport plans.
Numerical examples illustrate that the linear Gromov--Wasserstein distances, similarly as LOT,
can replace the expensive computation of pairwise  Gromov--Wasserstein distances
in applications \re{like shape classification}. 
\end{abstract}


\section{Introduction} \label{sec:intro}
\IEEEPARstart{R}{ecently,} a simplified version of optimal transport in Wasserstein
spaces, called linear optimal transport (LOT), 
was introduced by Wang et al.\ \cite{wang2013linear}.
\re{The} theoretical justifications \re{of LOT} can be found in the book of Ambrosio, Gigli and Savar\'e \cite{gradient_flows}.
From a geometric point of view, this approach just transfers measures 
from the geodesic Wasserstein space by the inverse exponential map
to the tangent space at some fixed reference measure that is 
assumed to be absolutely continuous (with respect to the Lebesgue measure). 
Then the LOT distance can be characterized by the optimal transport maps 
between the reference measure
and the considered measures.
This approach allows to work in the linear tangent space \re{rather than in the non-linear} Wasserstein space; 
so subsequent computations can utilize 
known methods from data science as for instance classification techniques.  
This is especially suited for
the approximate computation of pairwise distances for large databases
of images and signals.  Meanwhile LOT has been successfully applied
for several tasks in nuclear structure-based pathology
\cite{WOSSR2011}, parametric signal estimation \cite{RHNHR2020},
signal and image classification \cite{KTOR2016,PKKR2017}, modeling of
turbulences \cite{EN2020}, cancer detection
\cite{BKR2014,OTWCKBHR2014,TYKSR2015}, Alzheimer disease detection
\cite{PNAS21}, vehicle-type recognition \cite{GLDY2019} as well as for
de-multiplexing vortex modes in optical communications
\cite{PCNWDR2018}.  On the real line, LOT can further be written using
the cumulative density function of the random variables associated to
the involved measures.  This was used in combination with the Radon
transform under the name Radon-CDT \cite{KPR2016,PKKR2017}.  
We like to mention that
(inverse) exponential mappings were also used for the
iterative computation of Fr\'echet means, also known as barycenters, in Wasserstein
spaces in \cite{ABCM2016}.   
Furthermore, in \cite{ALR2021}, the determination of conditions that allow the
transformation of signals created by algebraic generative models into
convex sets by applying LOT has been addressed.  In
\cite{moosmueller2021linear}, the authors characterized settings in
which LOT embeds families of distributions into a space in which they
are linearly separable and provided conditions such that the LOT
distance between two measures is nearly isometric to the
Wasserstein distance. Finally, note that a linear version of the
Hellinger--Kantorovich distance is also available \cite{CCST2021}.

However, when dealing with reference measures that are not absolutely continuous, e.g.,
discrete measures, then optimal transport maps are in general not available such that
a generalized setting of LOT is needed. In this paper, we propose a generalized
LOT which relies on barycentric averaging maps of optimal transport plans
instead of optimal transport maps. 
For discrete measures, such an approach was also considered in \cite{wang2013linear}.
In this paper, we actually need this generalized LOT concept to motivate 
our framework of  (generalized) linear Gromov--Wasserstein distances.

Gromov--Wasserstein distances were first considered by M\'emoli in \cite{memoli2011gromov}
as a modification of Gromov--Hausdorff and Wasserstein distances.
A survey of the geometry of Gromov--Wasserstein spaces was given by Sturm in \cite{sturm2020space}.
Due to its invariance on isomorphism classes of so-called metric measure spaces,
the Gromov--Wasserstein distance is more suited for certain practical computations
like shape comparison and matching while retaining several desirable theoretical properties of its predecessors.
A combination with inverse problems has been considered in \cite{memoli2021distance}.
Further, a sliced version of the Gromov--Wasserstein distance has been discussed in \cite{sliced_gw,BHS21}. 
Recently, Gromov--Wasserstein distances were examined for Gaussian measures
in \cite{SDD2021}.

\IEEEpubidadjcol

In this paper, we introduce a linear variant of the Gromov--Wasserstein distance
that has the same advantages as LOT, namely the efficient computation of pairwise
distances in larger datasets, which can be subsequently coupled with standard
methods from image and signal processing.
Since the Brenier theorem that relates optimal transport maps with transport plans
in Wasserstein spaces is not available for the Gromov--Wasserstein setting,
we rely on optimal transport plans with respect to Gromov--Wasserstein distances
which always exist. 
Numerical examples illustrate \re{the excellent performance of the linear variant in shape classification tasks
and show that the distinctiveness remains comparable to 
the original Gromov--Wasserstein distance.}

\textbf{Outline of the paper:} 
In Section \ref{sec:ot}, we deal with linear optimal transport
and its generalization via barycentric projection maps of transport plans.
In Section \ref{linearGW}, we consider Gromov--Wasserstein distances.
We introduce the basic notation and properties that are quite technical, but we try to
keep things as simple as possible.
Then, following the definition of generalized LOT,  
we propose generalized linear Gromov--Wasserstein distances.
Section \ref{sec:numerical_examples} demonstrates how linear Gromov--Wasserstein distances
perform in several applications.
Finally, conclusions are drawn in Section \ref{sec:conclusions}.

\section{Linear Optimal Transport} \label{sec:ot}
In this section, we introduce a general version of LOT.
\re{
We will use the same underlying idea for the linear Gromov--Wasserstein distance.
}
\subsection{Optimal Transport}
By $L^2_\mu(\R^d,\R^d)$ we denote the space of (equivalence classes of) \re{measurable}
functions $T \colon \R^d \to \R^d$ fulfilling
\begin{equation*}
    \|T\|_{L^2_\mu} 
		\coloneqq 
		\Big( \int_{\R^d} \|T(x)\|^2 \dx \mu(x) \Big)^{\frac12} < \infty.
\end{equation*}
Let $\mathcal P(\R^d)$ be the space of probability
measures on the Borel $\sigma$-algebra $\mathcal B(\R^d)$, and 
$\mathcal P_2(\R^d)$ be the space of measures with finite second moments.
The \emph{push-forward measure}
$T_\# \mu$ of $\mu \in \mathcal P(\R^d)$ by a measurable map
$T: \R^d \to \mathcal \R^d$ is defined by
$T_\#\mu(B) \coloneqq \mu(T^{-1}(B))$ for all
$B \in \mathcal B(\R^d)$.
By $\| \cdot\|$ we denote the Euclidean norm on $\mathbb R^d$.
Together with  the \emph{Wasserstein distance} 
\begin{equation}\label{wasserstein}
  W(\mu,\nu) \coloneqq
   \min_{\pi \in \Pi(\mu,\nu)}
  \biggl( \int_{\R^d \times \R^d} \|x - y\|^2 \dx \pi(x,y) \biggr)^{\frac{1}{2}},
\end{equation}
where $\Pi(\mu,\nu)$ denotes the set of transport plans $\pi \in \mathcal P(\R^d \times \R^d)$
with marginals $\mu$ and $\nu$,
the space $\mathcal P_2(\R^d)$ becomes a metric space, known as (2-)\emph{Wasserstein space}. 
\re{
We denote the set of optimal transport plans, i.e. solutions to the minimization problem in \eqref{wasserstein}, by $\Pio(\mu,\nu)$.
}
For the more general definition of $p$-Wasserstein spaces, $p \in [1,\infty)$, 
see for instance \cite{villani2008optimal}.
The Wasserstein space is a \emph{geodesic space}
meaning that, for every $\mu,\nu \in \mathcal P_2(\R^d)$,
there exists a continuous curve $\gamma:[0,1]\to\mathcal P_2(\R^d)$ 
with $\gamma(0) = \mu$, $\gamma(1) = \nu$ and
\begin{equation} \label{geode}
W(\gamma(t),\gamma(s)) = |t-s| W(\gamma(0),\gamma(1)) 
\end{equation}
for all $t,s\in[0,1]$.
A continuous curve with property \eqref{geode} is called (constant speed) \emph{geodesic}.

If the measure $\mu$ is absolutely continuous, then, by the following theorem of
Brenier \cite{Brenier1987},
optimal transport plans in \eqref{wasserstein} are unique and 
can be characterized by transport maps.

\begin{theorem}[Brenier's Theorem]
  Let $\mu,\nu \in \mathcal P_2(\R^d)$, where $\mu$ is
  absolutely continuous. Then the minimization problem in \eqref{wasserstein} admits a unique solution $\pi^\nu_\mu$.
  Moreover, there exists a unique optimal transport map 
  $T^\nu_\mu \in L^2_\mu(\R^d,\R^d)$
  which solves
  $$
  \min_{T}
  \int_{\R^d} \|x - T(x)\|^2 \dx \mu(x) \quad \mathrm{subject \; to} \quad T_\# \mu = \nu.
  $$
  This optimal map is related to the optimal transport plan by
  $$\pi^\nu_\mu = (\id,T_\mu^\nu)_\# \mu.$$
\end{theorem}

The situation changes if $\mu$ is not absolutely continuous.
Then there still exists an optimal transport plan, but it may not be unique.
In contrast, the existence of an optimal transport map  is not guaranteed.
However, if there exists $T$ such that $\nu = T_\# \mu$ and $\pi \coloneqq (\id,T)_\# \mu$ is an optimal plan, then $T$ is an optimal map.
Conversely, if $T$ is an optimal map, then $\pi \coloneqq (\id,T)_\# \mu$ fulfills the marginal conditions, 
but must not be an optimal plan, as the example $\mu \coloneqq \frac14 \delta_0 + \frac34 \delta_1$ and 
$\nu \coloneqq \frac34 \delta_0 + \frac14 \delta_1$ shows.

\subsection{Linear Optimal Transport}
For discrete measures with a maximum of $n$ support points, 
the optimal transport amounts to solving a linear program that has worst-case complexity of \re{$n^3\log(n)$}.
Computing the pairwise Wasserstein distances of $N$ such measures results in $\binom{N}{2}$ optimal transport computations, which becomes numerically intractable for large $N$.
To speed up the numerical comparison, Wang et al.
\cite{wang2013linear} proposed LOT, 
which exploits the geometric structure of the Wasserstein space.
Following \cite[Eq~(8.5.1)]{gradient_flows}, 
the \emph{reduced tangent space (cone)} $\Tan^\red_\sigma \mathcal P_2(\mathbb R^d) \subset L^2_\mu(\mathbb R^d, \mathbb R^d)$ with base $\sigma \in \mathcal P_2(\mathbb R^d)$ is given by
\begin{align}
  &\Tan^\red_\sigma \mathcal P_2(\mathbb R^d) \notag
  \\
  &:=
    \overline{\bigl\{ r(T - \id) : (\id\re{,} T)_\# \sigma \in \Pio(\sigma, T_\# \sigma), r > 0 \bigr\}}^{L^2_\sigma}\!\!. 
\label{eq:red-tan}
\end{align} 
Note that the mapping $T$ in \re{\eqref{eq:red-tan}} is always an optimal transport map between $\sigma$ and $T_\# \sigma$.
If $\sigma$ is absolutely continuous, then 
the mapping 
$$F_\sigma: \mathcal P_2(\R^d) \to \Tan^\red_\sigma \mathcal P_2(\R^d), \quad \mu \mapsto T_\sigma^\mu-\id$$ 
is the inverse exponential map.
The key idea of LOT is to approximate $W(\mu,\nu)$ by the distance of the liftings to the tangent space, i.e.
\begin{equation}
  \label{eq:lot-ac}
  \LOT_\sigma(\mu, \nu) := \|F_\sigma(\mu) - F_\sigma(\nu)\|_{L_\sigma^2}
  = \|T_\sigma^\mu - T_\sigma^\nu  \|_{L_\sigma^2}.
\end{equation}
Then LOT is length preserving, i.e. $W(\mu,\sigma) = \LOT_\sigma(\mu, \sigma)$ and
 gives an upper bound of the Wasserstein distance
\[
W(\mu,\nu) \le \LOT_\sigma(\mu,\nu).
\]
For a fixed $\sigma \in \mathcal P_2(\mathbb R)$, the computation of all pairwise $\LOT_\sigma$ distances by \eqref{eq:lot-ac}
between $N$ measures requires only $N$ transport map computations.

One shortcoming of $\LOT_\sigma$ in \eqref{eq:lot-ac} is that the base measure $\sigma$ has to be absolutely continuous 
to ensure that the inverse exponential map to the reduced tangential space is well-defined for all measures in $\mathcal P_2(\mathbb R^d)$.
As a remedy, we replace the reduced tangent space by the geometric tangent space.
Given $\pi_\sigma^\mu \in \Pio(\sigma, \mu)$, the mapping
\begin{equation*} 
  t \mapsto
  \pi^{\sigma\to\mu}_t
  := ((1-t) P^1 + t P^2)_\# \pi_\sigma^\mu,
  \quad t \in [0,1],
\end{equation*} 
with the projections $P^1(s,x) \coloneqq s$ and $P^2(s,x) \coloneqq x$ defines a geodesic between $\sigma$ and $\mu$.
Moreover, every geodesic corresponds one-to-one to an optimal plan \cite[Thm~7.2.2]{gradient_flows}.
Henceforth, we identify each geodesic by its plan.
Let $G_\sigma$ denote the set of equivalence classes of all geodesics starting in $\sigma$, 
where two geodesics \re{$\pi^{\sigma\to\mu}_t$ and $\pi^{\sigma\to\nu}_t$} are equivalent if there exists an $\epsilon > 0$ such that \re{$\pi^{\sigma\to\mu}_t = \pi^{\sigma\to\nu}_t$} for $t \in [0,\epsilon]$.
The \emph{geometric tangent space} $\Tan^\geo_\sigma \mathcal P_2(\mathbb R^d)$ is the closure of $G_\sigma$ with respect to the metric
\begin{equation} \label{wsigma}
  W^2_\sigma(\pi_\sigma^\mu, \pi_\sigma^\nu)
  := \min_{\pi \in \re{\Gammas}(\pi_\sigma^\mu, \pi_\sigma^\nu)}
  \int_{\mathbb R^{3d}} |x - y|^2 \dx \pi(s,x,y),
\end{equation} 
where $\re{\Gammas}(\pi_\sigma^\mu, \pi_\sigma^\nu)$ consists of all 3-plans $\pi \in \mathcal P(\mathbb R^d \times \mathbb R^d \times \mathbb R^d)$ with $P^{12}_\# \pi = \pi_\sigma^\mu$ and $P^{13}_\# \pi = \pi_\sigma^\nu$, 
and where $P^{12}(s,x,y) \coloneqq (s,x)$ and $P^{13}(s,x,y) \coloneqq (s,y)$, cf.\ \cite[§~12.4]{gradient_flows}.  
Note that the plans $\pi \in \Gammas(\pi_\sigma^\mu, \pi_\sigma^\nu)$ also give rise to so-called generalized geodesics between $\mu$ and $\nu$, c.f.\ \cite[§~9.2]{gradient_flows}.

If $\sigma$ is not absolutely continuous, there may exist more than one geodesic between 
$\sigma$ and $\mu,\nu$, i.e.\ $\Pio(\sigma,\mu)$ and $\Pio(\sigma, \nu)$ are no singletons; 
so a proper extension of LOT to not absolutely continuous bases is
\begin{equation}
  \label{eq:lot-3p}
  \LOT_\sigma(\mu,\nu)
  := \inf_{\substack{\pi_\sigma^\mu \in \Pio(\sigma,\mu)\\ 
	\pi_\sigma^\nu \in \Pio(\sigma,\nu)}}
  W_\sigma(\pi_\sigma^\mu, \pi_\sigma^\nu).
\end{equation}  
\re{It can be verified} that LOT in \eqref{eq:lot-ac} and \eqref{eq:lot-3p} coincides for absolutely continuous $\sigma$.
In general $\LOT_\sigma$ is only a semi-metric, i.e., the triangular inequality is not fulfilled.
Taking the supremum instead of the infimum in \eqref{eq:lot-3p} would fix this issue.
Moreover, we have again $W(\mu,\nu) \le \LOT_\sigma(\mu, \nu)$.

\begin{remark}
Besides the geometric interpretation, we may interpret
$\mathrm{LOT}_\sigma$ as a constrained optimization of 
  \eqref{wasserstein}. More precisely, if we are given two plans
  $\pi^{\mu}_\sigma \in \Pio(\sigma, \mu)$ and
  $\pi^{\nu}_\sigma \in \Pio(\sigma, \nu)$
  \re{in \eqref{eq:lot-3p}}, 
	then the gluing lemma of Dudley
  \cite[Lem.~8.4]{ABS21} ensures the existence of
  $\pi_{\mathrm g} \in \mathcal P(\R^d \times \R^d \times \R^d)$
  such that 
	$P^{12}_\# \pi_{\mathrm g} = \pi_\sigma^{\mu}$ 
	and
  $P^{13}_\# \pi_{\mathrm g} = \pi_\sigma^{\nu}$.  
	The two plans $\pi_\sigma^{\mu}$ and
  $\pi_\sigma^{\nu}$ are glued together along the first axis.
	If the two   marginal plans are related to maps, i.e.\
  $\pi_\sigma^{\mu} = (\id, T_\sigma^{\mu})_\# \sigma$ and
  $\pi_\sigma^{\nu} = (\id, T_\sigma^{\nu})_\# \sigma$, 
	then the gluing is unique and given by
  $\pi_{\mathrm g} \coloneqq (\id, T_\sigma^{\mu}, T_\sigma^{\nu})_\# \sigma$.
  \re{Against this background,
  the marginal $P^{23}_\# \pi_{\mathrm g} \in \Pi(\mu,\nu)$  
  may be interpreted as transport from $\mu$ to $\nu$ via $\sigma$,
  and the optimization in \eqref{eq:lot-3p}
  is the constrained optimization of the Wasserstein distance \eqref{wasserstein} 
  restricted to the plans via $\sigma$.
  }
   \end{remark}
  
\subsection{Generalized Linear Optimal Transport}
Although $\LOT_\sigma$ in \eqref{eq:lot-3p} is also well defined for point reference measures, 
the numerical implementation requires the computation of an optimal 3-plan, which completely counteracts the intention behind LOT.
Instead we remain in the setting of transport maps by using
barycentric projection maps, 
which are based on the disintegration of transport plans \cite[Thm~5.3.1]{gradient_flows}.  
More precisely, given $\pi \in \mathcal P(\mathbb R^d \times \mathbb R^d)$ with $P^1_\# \pi = \sigma$, 
there exists a $\sigma$-almost everywhere uniquely defined family of measures $(\pi_s)_{s \in \R^d} \subset \mathcal P(\mathbb R^d)$ such that
\begin{equation*}
  \int_{\mathbb R^{2d}} f(s,x) \dx \pi(s,x)
  = \int_{\mathbb R^d} \int_{\mathbb R^d} f(s,x) \dx \pi_s(x) \dx \sigma(s).
\end{equation*}
\re{for all measurable functions $f:\R^{2d} \to [0,\infty)$.}
The \emph{barycentric projection map} $\mathcal T_{\pi}: \R^d \to \R^d$
of $\pi \in \mathcal P(\mathbb R^d \times \mathbb R^d)$ with first marginal $\sigma$
is defined for $\sigma$-almost every $s \in \mathbb R^d$ by 
\begin{equation} \label{eq:baryproj}
  \mathcal T_{\pi}(s) \coloneqq  \int_{\mathbb R^d} x \dx \pi_s(x) = \argmin_{x' \in \mathbb R^d} \int_{\mathbb R^d} \|x-x'\|^2 \dx \pi_s(x)
\end{equation} 
provided that $\pi_s$ has finite second moments $\sigma$-a.e., see \cite[p.
126]{gradient_flows}.

\begin{example}
\re{Let $\delta_x \in \p(\R^d)$ and $\delta_{(s,x)} \in \p(\R^d \times \R^d)$ denote the Dirac measure at $x \in \R^d$ and $(s,x) \in \R^d \times \R^d$ respectively}. For the discrete probability measures 
\[
\sigma = \sum_{i=1}^n \sigma_i \delta_{s_i} \in \mathcal P(\mathbb R^d)
\quad\re{\text{and}\quad
\mu = \sum_{j=1}^m \mu_j \delta_{x_j} \in \mathcal P(\mathbb R^d)}
\]
and \re{the transport plan}
\[
\pi = \sum_{i=1}^n \sum_{j=1}^m \pi_{i,j} \delta_{\re{(}s_i,x_j\re{)}} 
\in \re{\Pi(\sigma,\mu)}
\]
with \re{$\sum_{j=1}^m \pi_{i,j} = \sigma_i$ and $\sum_{i=1}^n \pi_{i,j} = \mu_j$}, 
the barycentric projection reads as
\[
\T_{\pi}(s_i) = \frac{1}{\sigma_i}\sum_{j=1}^m \pi_{i,j} x_j, \quad i = 1,\dotsc,n.
\]
Such maps are also used in \cite{wang2013linear}.
\end{example}

By the following proposition, the barycentric projection map \eqref{eq:baryproj} of an optimal transport plan $\pi_\sigma^\mu$ 
is always an optimal transport map $T_\sigma^{\tilde \mu}$ between
$\sigma$ and $\tilde \mu = (\mathcal T_{\pi_\sigma^\mu})_\# \sigma$.

\begin{proposition} \label{prop_1} For each $\pi_\sigma^\mu \in \Pio(\sigma,\mu)$, the barycentric projection map 
$\mathcal T_{\pi_\sigma^\mu}: \mathbb R^d \to \mathbb R^d$ in \eqref{eq:baryproj} 
defines an optimal transport map from $\sigma$ to the measure $\tilde \mu \coloneqq (\mathcal T_{\pi_\sigma^\mu})_\# \sigma$, i.e., 
$$\mathcal T_{\pi_\sigma^\mu} = T^{\tilde \mu}_\sigma. $$
\end{proposition}

Although the statement may be implicitly derived from \cite[§~12.4]{gradient_flows}, we give a direct proof in the appendix.
On the basis of the barycentric projection, 
we propose to extend the LOT formulation in \eqref{eq:lot-ac} by considering \emph{generalized LOT (gLOT)}
\begin{equation}
    \label{eq:w-lot}
  \gLOT_\sigma (\mu, \nu)
  :=
  \inf_{\substack{\pi_\sigma^\mu \in \Pio(\sigma,\mu) \\ \pi_\sigma^\nu \in \Pio(\sigma, \mu)}}
  \| \mathcal T_{\pi_\sigma^\mu} - \mathcal T_{\pi_\sigma^\nu} \|_{L^2_\sigma}.
\end{equation}
If $\pi_\sigma^\mu = (\id,T_\sigma^\mu)_\# \sigma$, then $\T_{\pi_\sigma^\mu} = T_\sigma^\mu$,
so that gLOT coincides with LOT in particular for absolutely continuous bases.
In the numeric\re{al} implementation of gLOT, the minimization over $\Pio(\sigma,\mu)$ and $\Pio(\sigma,\nu)$ \re{in \eqref{eq:w-lot} can be omitted}, i.e.\ we use fixed transport plans $\pi_\sigma^\mu$ and $\pi_\sigma^\nu$ instead.

\begin{remark}
gLOT has actually a geometric interpretation.  
The barycentric projection $\pi \to \mathcal T_\pi$ defines a map from the geometric tangent space to the reduced tangent space by
\begin{equation*}
  \pi_\sigma^\mu \in \Tan_\sigma^\geo \mathcal P_2(\mathbb R^d)
  \quad
  \mapsto
  \quad
  (\mathcal T_{\pi_\sigma^\mu} - \id) \in \Tan_\sigma^\red \mathcal P_2(\mathbb R^d), 
\end{equation*} 
see Proposition \ref{prop_1}.
From this point of view, gLOT takes two geodesics corresponding to the optimal plans $\pi_\sigma^\mu$ and $\pi_\sigma^\nu$ from the geometric tangent space, maps them to the reduced tangent space, and computes the distance there.
In this way, we overcome the issue that the inverse exponential map may not be defined for the whole $\mathcal P_2(\mathbb R^d)$, 
\re{which prevents the application of \eqref{eq:lot-ac} in the discrete setting,
and the issue of the costly computation of \eqref{eq:lot-3p}.}
\end{remark}

\section{Linear Gromov--Wasserstein Distance}  \label{linearGW}

\re{In certain applications like shape matching, 
the Wasserstein distance is unfavourable since it varies under
isometric transformations such as
translations and rotations of the considered measures. 
For this reason,}  
Mémoli
\cite{memoli2011gromov} introduced an optimal-transport-like distance,
where the aim was to match measures according to pairwise distance
perturbations.  To this end, we need the definition of a \emph{metric
  measure space (mm-space)}, which is a triple $\XX = (X, d_{X}, \mu)$,
where
\begin{itemize}
\item[i)] $(X, d_{X})$ is a compact metric space,
\item[ii)] $\mu$ is a Borel probability measure on $X$ with full
  support.
\end{itemize}

\subsection{Gromov--Wasserstein Distance}
For two mm-spaces  $\XX = (X,d_{X},\mu)$ and $\YY =
(Y,d_{Y}, \nu)$, 
the \emph{Gromov--Wasserstein (GW) distance} is defined by
\begin{align}
&\GW(\XX,\YY) 
\coloneqq\hspace*{-10pt} 
\inf_{\pi \in \Pi(\mu,\nu)} 
\biggl(  \int_{(X \times Y)^2} \hspace{-25pt} \lvert d_{X}(x,x') - 
d_{Y}(y,y') \rvert^2 \\[-8pt]
  &\hspace{110pt} \times\dx \pi(x,y) \dx \pi(x',y') \biggr)^{\tfrac12}.
    \label{eq:GW}
\end{align}
Here $\pi \in \Pi(\mu,\nu)$ means that $\pi \in \mathcal P(X\times Y)$ 
has marginals $\mu$ and $\nu$.
Further, we denote by $\Pio(\XX,\YY)$ the set of optimal GW plans in \eqref{eq:GW}.
In the literature, the above quantity is also called
2-Gromov--Wasserstein distance, and analogous definitions for
$p \in [1,\infty)$ as well as further generalizations are possible.
For an overview, we refer also to \cite{Vayer2020}.
Due to the
Weierstraß theorem, a minimizer in \eqref{eq:GW} always exists
\cite[Cor~10.1]{memoli2011gromov}.  Two mm-spaces
$\XX = (X,d_{X},\mu)$ and
$\YY = (Y,d_{Y},\nu)$ are called \emph{isomorphic} if and
only if there exists a (bijective) isometry $\psi: X \to Y$
such that $\psi_\# \mu = \nu$.  
We denote the corresponding equivalence classes by $\llbracket \cdot \rrbracket$.
The GW distance defines a metric on these equivalence classes  \cite[Thm~5.1]{memoli2011gromov}.
The resulting (incomplete) metric space is  called the \emph{Gromov--Wasserstein space}.  In
particular, the GW distance is invariant under
translation and rotation of the mm-space.  

\re{Up to now, there does not exist a general GW analogue to Brenier's Theorem, which would ensure the existence of optimal plans that are realized by transport maps under certain regularity assumptions.
A comprehensive overview on this specific subject is given in \cite[Rem~3.3]{memoli2021distance}}. In \cite{sturm2020space}, Sturm has shown that in the Euclidean setting 
optimal GW plans between rotationally invariant probability spaces are realized by optimal transport maps.

Due to its invariance properties and independence of the ambient spaces, 
the GW metric provides a valuable tool for data science, 
shape analysis, and object classification. 
However, its exact computation is NP-hard. 
Even its approximation is computationally challenging and requires, if
a gradient descent algorithm is used,
$O(n^3 \log(n))$ arithmetic operations, where $n$ is the cardinality of the underlying mm-spaces \cite{PCS2016}. 
For improvements in the setting of sparse graphs, see \cite{XLZD2019}.
Hence its use for comparing a larger number of mm-spaces is limited,
which motivates the following considerations.

\subsection{Linear Gromov--Wasserstein Distance}
We consider the (equivalence classes of) mm-spaces
$\SS = \llbracket S,d_S,\sigma \rrbracket$, $\XX = \llbracket X,d_X,\mu \rrbracket$, 
and $\YY = \llbracket Y,d_Y,\nu \rrbracket$.
In contrast to the above definition from Mémoli
\cite{memoli2011gromov}, we  allow that the measures $\sigma$,
$\mu$, and $\nu$ may  not have full support.
Similarly to the Wasserstein setting, the Gromov--Wasserstein space is geodesic.
\re{The construction of the tangent space is, however, more technical.} 
We follow the lines of Sturm in  \cite{sturm2020space}.
Each geodesic from $\SS$ to $\XX$ has the form
\begin{equation}
  \label{eq:gw-geo}
  t \mapsto \pi_t^{\SS\to\XX}
  := \llbracket S \times X, (1-t) d_S + t d_X, \pi \rrbracket,
  \quad
  t \in [0,1],
\end{equation}
where $\pi \in \Pio(\SS, \XX)$, 
and 
where $d_S$ acts on the $S$ components and $d_X$ on the $X$ components of $(S \times X)^2$, 
respectively.
Conversely, every optimal plan defines a geodesic.
Note that $\pi_0^{\SS\to\XX}$ and $\pi_1^{\SS\to\XX}$ are isomorphic to $\SS$ and $\XX$ by  
$P^1(s,x)\coloneqq s$ and $P^2(s,x) \coloneqq x$, respectively. 

In order to introduce tangent spaces and to derive their explicit representations, 
the GW space is embedded into the more regular space of gauged measure spaces.
A \emph{gauged measure space (gm-space)} is as before a triple $\re{\FS} :=(S,k_S,\sigma)$,
where the distance is replaced by a so-called \emph{gauge function} 
$k_S$ in $L^2_{\re{\sym}}(S \times S, \sigma \otimes \sigma)$, 
\re{which} consists 
of all \emph{symmetric}, square-integrable functions with respect to $\sigma \otimes \sigma$.
Here, $S$ can be a Polish space.
\re{
Note that gm-spaces are more general than mm-spaces as gauge functions include,
for instance, pseudometrics (which are not definite) and semimetrics (which do not admit the triangle inequality) on compact spaces.
Clearly, every mm-space is a gm-space.
}
The extension of the GW distance to the gm-spaces 
$\re{\FX} = (X, k_X, \mu)$ and 
$\re{\FY} = (Y, k_Y, \nu)$ is given by
\begin{align}
  &\GW(\re{\FX},\re{\FY}) = \inf_{\pi \in \Pi(\mu,\nu)}
    \biggl({\int\limits_{(X \times Y)^2}}
    \lvert k_X(x,x') - k_Y(y,y') \rvert^2
    \notag \\[-10pt]
  \label{eq:gw-gm}
  &\hspace{120pt} \times \dx \pi(x,y) \dx \pi(x',y')
    \biggr)^{\frac12} \!\!\!.
\end{align}
A minimizing coupling always exists \cite[Thm~5.8]{sturm2020space}.
The set of all plans minimizing the integral in \eqref{eq:gw-gm} with respect to 
$\re{\FX}$ and $\re{\FY}$ 
is henceforth denoted by $\Pio(\re{\FX},\re{\FY})$.
Two gauged measure spaces $\re{\FX}$ and $\re{\FY}$ are called \emph{homomorphic} 
if $\GW(\re{\FX},\re{\FY}) = 0$.
The space $\mathfrak G$ of 
homomorphic equivalent classes---again denoted by $\llbracket \cdot \rrbracket$---equipped 
with the GW distance \eqref{eq:gw-gm} is complete and geodesic.
To simplify notation, we denote such equivalence classes again  by $\re{\FX}$.
The geodesics from $\re{\FS}$ to $\re{\FX}$ have the form
\begin{equation}
  \label{eq:gm-geo}
  t \mapsto \pi_t^{\re{\FS}\to\re{\FX}}
  := \llbracket S \times X, (1-t) k_S + t k_X, \pi \rrbracket,
  \quad
  t \in [0,1],
\end{equation}
where $\pi \in \Pio(\re{\FS},\re{\FX})$. 
Conversely, each $\pi \in \Pio(\re{\FS},\re{\FX})$ defines a geodesic.

Formally, the \emph{tangent space} $\Tan_{\re{\FS}} \mathfrak G$ with base 
$\re{\FS} \in \mathfrak G$ is defined as
\begin{equation*}
  \Tan_{\re{\FS}} \mathfrak G 
	\coloneqq 
	\Big( \bigcup_{\llbracket  S, k_S ,\sigma \rrbracket = \re{\FS}} 
	L_{\re{\sym}}^2(S \times  S, \sigma \otimes \sigma) \Big) {\Bigm/}{\sim},
\end{equation*}
where the union is taken over all gm-spaces $(S, k_S, \sigma)$ in the equivalence class $\re{\FS}$
and two functions 
$g \in L_{\re{\sym}}^2(S \times S,\sigma \otimes \sigma)$ 
and 
$g' \in L_{\re{\sym}}^2( S' \times  S',\sigma' \otimes \sigma')$ 
defined on the representatives 
$(S, k_S, \sigma)$ and $(S', k'_S, \sigma')$ of $\re{\FS}$
are equivalent, 
if there exists  $\pi \in \Pio \left((S ,k_{S},\sigma), (S' ,k_{S'},\sigma' ) \right)$ 
such that
\begin{equation*}
  g(s_1,s_2) = g'(s_1',s_2')
\end{equation*}
almost everywhere with respect to $\pi(s_1,s_1')\otimes\pi(s_2,s_2')$.
Note that each tangent 
$g \in \Tan_{\re{\FS}} \mathfrak G$ 
is implicitly associated with its representative $(S, k_S, \sigma)$.
A  (cone) \emph{metric on} $\Tan_{\re{\FS}}   \mathfrak G$ is given by
\begin{align}
  \notag &\GW_{\re{\FS}}  (g,h)
  \\
  \label{eq:gw-tan}
         &\coloneqq
           \inf \bigl\{\|g-h\|_{L^2((S \times S')^2,\pi \otimes \pi)} :
           \pi \in \Pio(\TT_g,\TT_h) \bigr\},
\end{align}
where $\TT_g$ and $\TT_h$ denote the representatives associated with $g$ and $h$.
Given $g \in \Tan_{\re{\FS}} \mathfrak G$ defined on the representative 
$(S, k_S, \sigma)$ of the equivalence class $\re{\FS}$, 
the \emph{exponential map} $E_{\re{\FS}} : \Tan_{\re{\FS}} \mathfrak G \to \mathfrak G$ is defined by
\begin{equation*}
  E_{\re{\FS}}(g) = \llbracket  S, k_S + g, \sigma \rrbracket.
\end{equation*}
As a consequence, every geodesic in \eqref{eq:gm-geo} may be written as
\begin{equation*}
  \pi_t^{\re{\FS}\to\re{\FX}} = E_{\re{\FS}}(t h)
  \quad\text{with}\quad
  h := k_X - k_S,
\end{equation*}
where $h$ is defined on the representative  $( S \times  X, k_S, \pi)$ with $\pi \in \Pio(\re{\FS},\re{\FX})$.
Note that two geodesics which coincide for all $t \in [0,\epsilon]$ for some $\epsilon > 0$ correspond to the same tangent; 
so the tangent space embrace all geodesics starting in $\re{\FS}$.
Associating any geodesic $\pi_t^{\re{\FS}\to\re{\FX}}$ 
with its optimal plan 
$\pi_{\re{\FS}}^{\re{\FX}} \in \Pio(\re{\FS},\re{\FX})$, 
we define $F_{\re{\FS}}: \mathfrak G \to \Tan_{\re{\FS}} \mathfrak G$ by
\begin{equation*}
  F_{\re{\FS}} (\pi_{\re{\FS}}^{\re{\FX}} ) = k_X - k_S
  \quad
   (\text{acting on} \;
  (S \times X, k_S, \pi_{\re{\FS}}^{\re{\FX}})).
\end{equation*}
For the geodesics \eqref{eq:gw-geo} between mm-spaces, we especially have
\begin{equation}
  \label{eq:lift-geo}
  F_{\SS}(\pi_{\SS}^{\XX}) = d_X - d_S
  \quad
  (\text{acting on} \;
  (S \times X, d_S, \pi_{\SS}^{\XX})).
\end{equation}
Against this background, we define the distance between two geodesics 
$\pi_{\SS}^{\XX}$ 
and 
$\pi_{\SS}^{\YY}$ as
\begin{equation}\label{dist_final}
  \GW_{\SS}(\pi_{\SS}^{\XX}, \pi_{\SS}^{\YY})
  \coloneqq
  \GW_{\SS}( F_{\SS}(\pi_{\SS}^{\XX}), F_\SS(\pi_{\SS}^{\YY})).
\end{equation}
Then we have the following relation whose proof is given in the appendix.

\begin{proposition}   \label{thm:exp-gw-geo}
  Consider the mm-spaces $\SS = \llbracket S,d_S,\sigma \rrbracket$, $\XX = \llbracket X,d_X,\mu \rrbracket$, $\YY = \llbracket Y,d_Y,\nu \rrbracket$.  The distance \eqref{dist_final} between the geodesics related to  $\pi_\SS^\XX \in \Pio(\SS,\XX)$ and $\pi_\SS^\YY \in \Pio(\SS,\YY)$ is given by
  \begin{align}\label{eq:GW_S}
    &\GW_\SS^2(\pi_\SS^\XX, \pi_\SS^\YY)
      = \inf_{\pi \in \re{\Gamma}_\SS(\pi_\SS^\XX, \pi_\SS^\YY)}
      \int_{(S \times X \times Y)^2}
      \hspace{-35pt}
      \lvert d_X(x,x') - d_Y(y,y') \rvert^2
    \nonumber\\
    &\hspace{105pt}
      \times \dx \pi(s,x,y) \dx \pi(s',x',y'),
  \end{align}
  where $\re{\Gamma}_\SS(\pi_\SS^\XX, \pi_\SS^\YY)$ consists of all 3-plans $\pi \in \mathcal P(S \times X \times Y)$ 
	with $P^{12}_\# \pi = \pi_\SS^\XX$ and $P^{13}_\# \pi = \pi_\SS^\YY$.
\end{proposition}

The minimization over the 3-plans in \re{\eqref{eq:GW_S}}
\re{is analogous to the minimzation over the} 3-plans in the definition of $W_\sigma$ in \eqref{wsigma}.
In the spirit of LOT, we now propose to approximate the GW distance $\GW(\XX, \YY)$ 
by lifting $\XX$ and $\YY$ to $\Tan_\SS \mathfrak G$ via geodesics and using the metric on the tangent space. 
More precisely, we define the \emph{linear Gromov--Wasserstein distance}  by
\begin{equation}\label{eq:GW_sigma}
  \LGW_{\SS}(\XX,\YY)
  \coloneqq \inf_{\substack{\pi_\SS^\XX \in \Pio(\SS,\XX)\\ \pi_\SS^\YY \in \Pio(\SS,\YY)}}
  \GW_\SS(\pi_\SS^\XX, \pi_\SS^\YY).
\end{equation}
In comparison with \eqref{eq:lot-3p}, we can consider LGW as an analogue to LOT
in \re{ the GW space}.
\re{Further, we have the following lower and upper bound, whose proof is given in the appendix.}

\re{
\begin{lemma}
\label{lem:LGW_estimates}
Let $\SS,\XX,\YY$ be mm-spaces. Then LGW is bounded above and below by
\begin{equation}
    \label{eq:LGW_estimates}
    \GW(\XX,\YY) \leq \LGW_\SS(\XX,\YY) \leq \GW(\SS,\XX) + \GW(\SS,\YY).
\end{equation}
\end{lemma}
}

\re{
\begin{remark}\label{rem:ref_space}
    The quality of the approximation of GW by LGW
    crucially depends on the chosen reference space $\SS$.
    In the sense of Lemma~\ref{lem:LGW_estimates},
    a suitable reference $\SS$ should ensure a small right-hand side in \eqref{eq:LGW_estimates}.
    For the approximation of the pairwise GW distances of
    several mm-spaces $\XX_1,\dotsc,\XX_N$,
    an appropriate reference $\SS$ should thus be equally close to all $\XX_k$.
    Since the minimization of $\sum_{k=1}^N \GW(\SS,\XX_k)$ over all mm-spaces $\SS$ 
    is intractable,
    a possible alternative would be a Gromov--Wasserstein barycenter
    minimizing $\sum_{k=1}^N \GW^2({\SS},\XX_k)$,
    which is discussed in more detail during the numerical experiments in Section~\ref{sec:numerical_examples}.
\end{remark}
}

\subsection{Generalized Linear Gromov--Wasserstein Distance}
Assume for the moment that $\pi_{\SS}^\XX\in \Pio(\SS,\XX)$ 
and $\pi_\SS^\YY \in \Pio(\SS,\YY)$ are unique and induced by optimal maps $T_\SS^\XX$ and $T_\SS^\YY$.
In this situation, $\re{\Gamma}_\SS(\pi_\SS^\XX,\pi_\SS^\YY)$ becomes the singleton   $(\id, T_\SS^\XX, T_\SS^\YY)_\# \sigma$, and we obtain
\begin{align}
  &\GW_{\SS}(\mathbb X, \mathbb Y) \notag
  \\
  &= \bigl\| d_{X}(T_\SS^\XX(\cdot_1), T_\SS^\XX(\cdot_2))
    - d_{Y}(T_\SS^\YY(\cdot_1), T_\SS^\YY(\cdot_2))
    \bigr\|_{L^2_{\sigma \otimes \sigma}}\!\!,\quad
    \label{eq:GWS}
\end{align}
where  $\cdot_1$ and $\cdot_2$
are the first and second argument with respect to
$S \times S$.

Similarly to LOT, LGW does not alleviate the computational costs of calculating pairwise $\GW$ distances.
For this reason, we recommend to 
use the barycentric projection mapping to transform $\pi_\SS^\XX$ and $\pi_\SS^\YY$ into maps 
$\mathcal T_{\pi_\SS^\XX}$ and $\mathcal T_{\pi_\SS^\YY}$.
Since the metric spaces may be more general 
than the measure spaces considered in Section~\ref{sec:ot}, we introduce the \emph{generalized barycentric projection}
\begin{equation}
  \label{eq:gen-bary-proj}
  \mathcal T_{\pi_\SS^{\XX}}(s)
  \coloneqq 
	\argmin_{x' \in  X} \int_{X}
  {d_{X}^2(x',x)} \, \dx \pi_{\SS,s}^\XX (x),
\end{equation}
where $\pi_{\SS,s}^\XX$ is the disintegration of the chosen
$\pi_\SS^\XX$. Based on the Weierstraß theorem, the minimum is attained.
In the special case that $X \subset \R^d$ is convex and
$d_X(x_1, x_2) = \|x_1 - x_2\|$, the generalized barycentric projection
coincides with \eqref{eq:baryproj}.

Analogously to gLOT \re{in \eqref{eq:w-lot} and based on $\GW_\SS$ in \eqref{eq:GWS}}, we define \emph{generalized LGW (gLGW)} by
\begin{align}    \label{eq:gLGW}
  &\gLGW_\SS(\XX, \YY)
  \\
  &\coloneqq
    \inf_{\substack{\pi_\SS^\XX \in \Pio(\SS,\XX)\\ \pi_\SS^\YY \in \Pio(\SS,\YY)}}
  \bigl\| d_{X}(\mathcal T_{\pi_\SS^\XX}(\cdot_1), \mathcal T_{\pi_\SS^\XX}(\cdot_2))
  \\[-15pt]
  &\hspace{110pt}
    - d_{Y}(\mathcal T_{\pi_\SS^\YY}(\cdot_1), \mathcal T_{\pi_\SS^\YY}(\cdot_2))
    \bigr\|_{L^2_{\sigma \otimes \sigma}}. \nonumber
\end{align}
For numerical computations, we again propose to use fixed optimal plans instead of minimizing over $\Pio(\SS,\XX)$ and $\Pio(\SS,\YY)$.

\begin{remark}
  \re{As mentioned in the introduction}, under the conditions of the Brenier theorem, the
  OT and LOT distances coincide in the one-dimensional setting.
  The linear GW distance differs in general from the GW distance
  also in one dimension.  We verified this by computing the
  corresponding GW and LGW distances for
  \begin{align*}
    S = \{0,1,2,3,6\}, \quad
    X = \{0,1,2,5,7\}, \quad
    Y = \{0,2,3,6,7\}
  \end{align*}
  the absolute value distances and the corresponding discrete measures
  with weights $\frac15$.  For this specific instance, we obtain
  $\GW(\XX, \YY) \approx 0.69$ and $\gLGW(\XX,\YY) \approx 1.13$.
\end{remark}

\section{Numerical Examples}\label{sec:numerical_examples}

All numerical experiments%
\footnote{\re{The source code is publicly available at \url{https://github.com/Gorgotha/LGW}.}}
in this section \re{have been performed on an off-the-shelf MacBook Pro (Apple M1 chip, 8GB RAM) and} have been implemented in
\emph{Python~3}, where we mainly rely on the packages \emph{Python
  Optimal Transport (POT)} \cite{POT-toolbox}, \emph{scikit-learn}
\cite{scikit-learn}, and \emph{NetworkX} \cite{networkx}.  POT
contains a Gromov--Wasserstein module allowing the numerical
computation of the GW distance \eqref{eq:GW}, a
corresponding optimal plan, and GW barycenters for
discrete mm-spaces, where the measure space consists of finitely many
points, and where the measure thus becomes a point measure. \re{A}
GW barycenter $\SS$ between the discrete mm-spaces
$\XX_k$ for $k=1,\dots,K$ is defined via
\begin{equation}
    \label{eq:bary-num}
  \SS \re{\in} \argmin_{\tilde{\SS}} \sum_{k=1}^K \GW^2(\tilde{\SS}, \XX_k).
\end{equation}
The minimization here goes over the set of all discrete mm-spaces with a certain number of points.  
\re{The corresponding POT method additionally presets the weights in $\SS$ 
and only minimizes over the metric $d_\SS$.}
To visualize the computed
GW barycenters and the computed pairwise gLGW
distances, we use the scikit-learn implementation of multi-dimensional
scaling (MDS) from \cite{scikit-learn}, which allows to embed a series
of points with given distances into $\R^2$ such that the
distances are approximately preserved.

\begin{figure}[t]
    \centering
    \includegraphics[width = 0.5\linewidth]{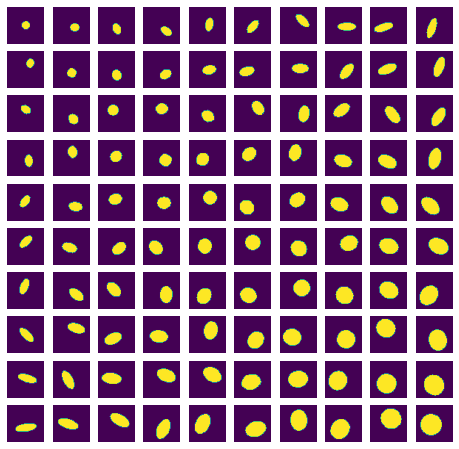}
    \caption{Elliptical disks used for the first numerical
      experiment.  The opposing images with respect to the diagonal
      form isometrical pairs resulting from rotations and shifts.}
    \label{fig:ells_data}
\end{figure}

\begin{table*}
\centering
\caption{\re{Quantitative comparison between the computed gLGW distances.
The first row shows the employed reference spaces $\SS_1,\dots, \SS_9$,
which include uniform distributions on different shapes as well as 
non-uniform distributions on the square $\SS_3$ and the two circles $\SS_8$.
The distribution is indicated by the color of the pixels.
For each reference, the computation time to compute all pairwise distances, 
the Mean Relative Error (MRE), and the Pearson Correlation Coefficient (PCC)
as well as the number of non-zero points in the reference are recorded.}}
\label{tab:1}
\re{
\begin{tabular}{lcccccccccc}
    \toprule
    reference & GW & $\SS_1$ & $\SS_2$ & $\SS_3$ & $\SS_4$ & $\SS_5$ & $\SS_6$ & $\SS_7$ & $\SS_8$ & $\SS_9$ \\
    && \includegraphics[width=30pt, trim=115pt 35pt 100pt 35pt, clip=true]{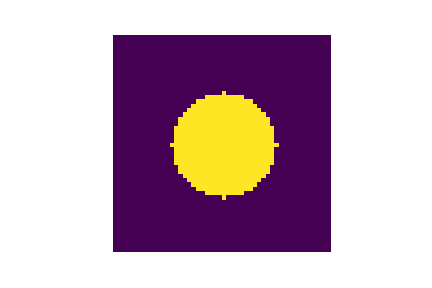}
    & \includegraphics[width=30pt, trim=115pt 35pt 100pt 35pt, clip=true]{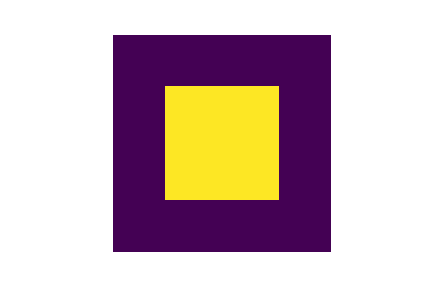}
    & \includegraphics[width=30pt, trim=115pt 35pt 100pt 35pt, clip=true]{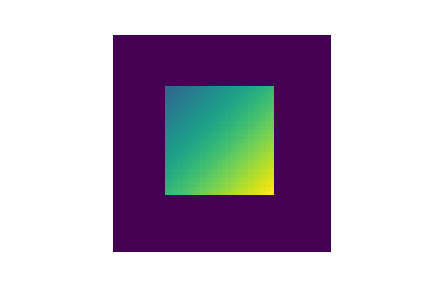}
    & \includegraphics[width=30pt, trim=115pt 35pt 100pt 35pt, clip=true]{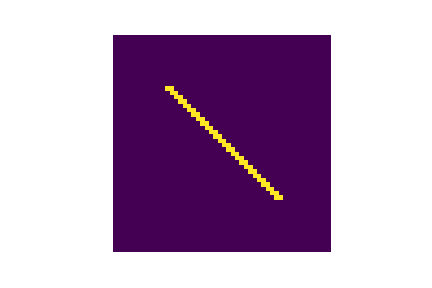}
    & \includegraphics[width=30pt, trim=115pt 35pt 100pt 35pt, clip=true]{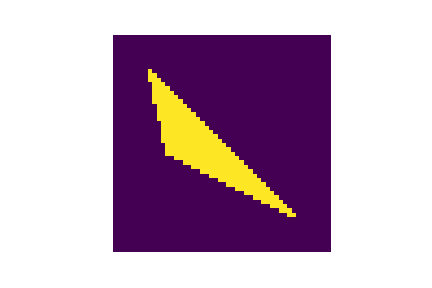}
    & \includegraphics[width=30pt, trim=115pt 35pt 100pt 35pt, clip=true]{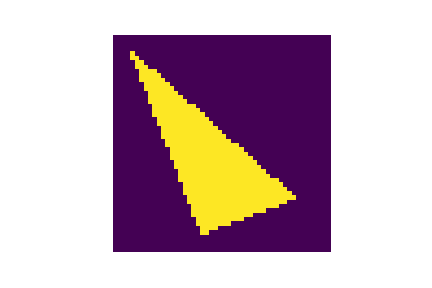}
    & \includegraphics[width=30pt, trim=115pt 35pt 100pt 35pt, clip=true]{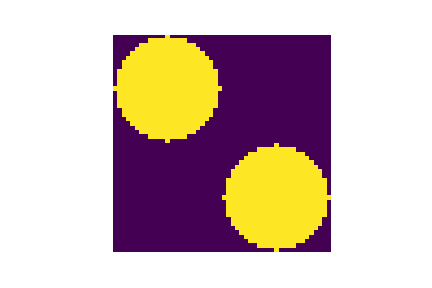}
    & \includegraphics[width=30pt, trim=115pt 35pt 100pt 35pt, clip=true]{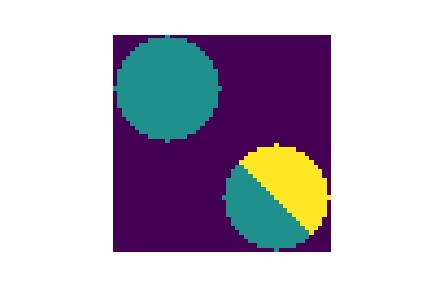}
    & \includegraphics[width=30pt, trim=115pt 35pt 100pt 35pt, clip=true]{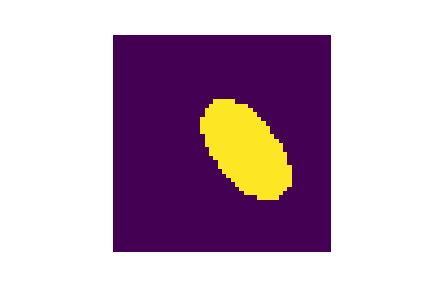}
    \\
    \midrule
    time & 159.38 min & 2.78 min & 15.92 min & 11.77 min & 0.13 min & 1.54 min & 5.45 min & 7.68 min & 8.28 min & 2.06 min
    \\
    MRE & --- & 0.336 & 0.325 & 0.312 & 0.158 & 0.038 & 0.030 & 0.017 & 0.016 & 0.019
    \\
    PCC & --- & 0.891 & 0.876 & 0.887 & 0.986 & 0.999 & 0.999 & 0.999 & 0.999 & 0.999
    \\
    points & --- & 441 & 676 & 625 & 52 & 289 & 545 & 882 & 882 & 317
    \\
    \bottomrule
\end{tabular}
}
\end{table*}

\begin{figure*}
    \centering
    \includegraphics[width = \linewidth]{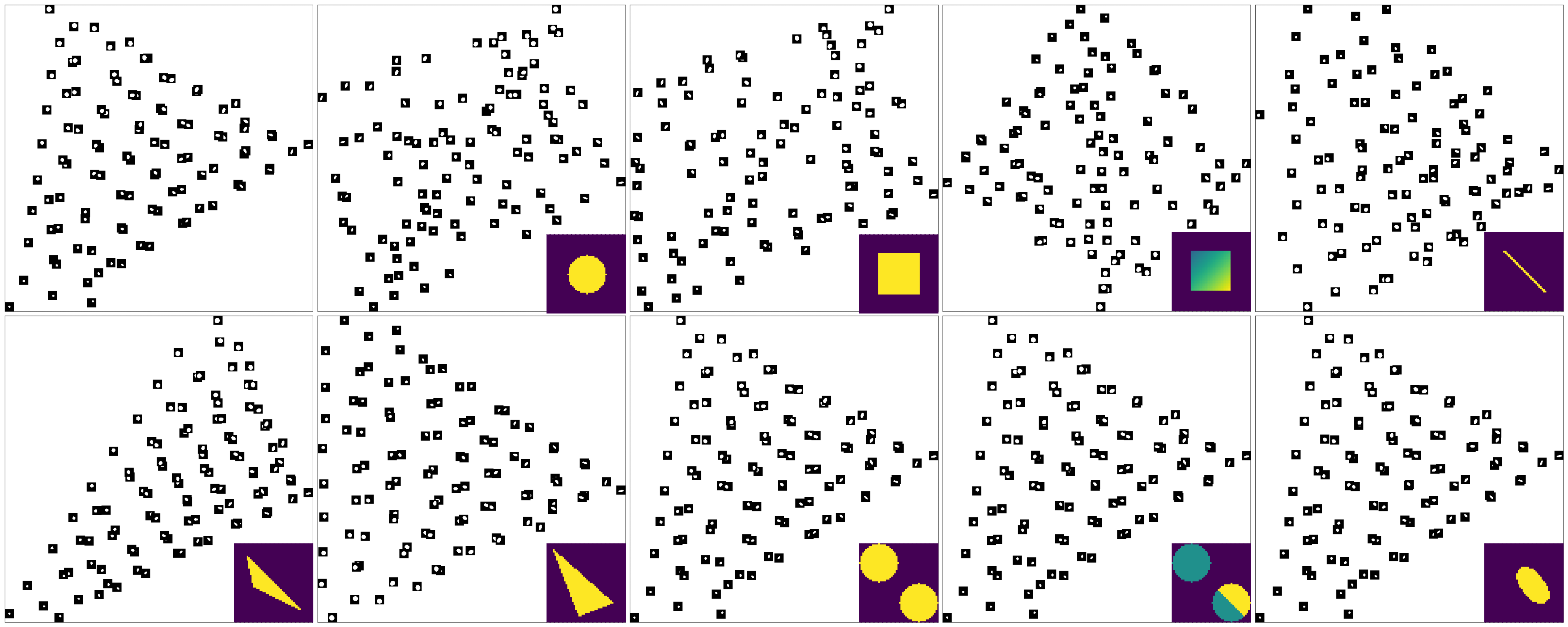}
    \caption{\re{MDS embdding of the computed pairwise GW distances (top left)
    and the gLGW distances based on the references in Table~\ref{tab:1}.
      }}
    \label{fig:ells_dists}
\end{figure*}

\subsection{Gromov--Wasserstein of elliptical disks}

For our first example, we apply $\gLGW$ to a toy problem, where we
want to compute the GW distance between a series of elliptical disks,
see Figure~\ref{fig:ells_data}.  Each image here consists of
$50 \times 50$ equispaced pixels in $[0, 50]^2$.  For the
numerical simulations, we interpret these images as discrete mm-spaces
$\XX_1, \dots, \XX_{100}$. For this  we set
$\XX_k \coloneqq([0,50]^2, d_{\mathrm E}, \mu_k)$, where $d_{\mathrm E}$ is the Euclidean distance
and $\mu_k$
corresponds the uniform distribution on the position of the white
pixels.  Notice that, except for the diagonal, all elliptical disks in
Figure~\ref{fig:ells_data} occur in isometrical pairs (up to
discretization errors).  Since the GW distance is invariant under
isometries, this should be reflected in the computed $\GW$ and $\gLGW$
distances.

For comparison, we first compute all pairwise GW
distances, where we use the optimal Wasserstein coupling as starting
value for the corresponding POT algorithm in the GW distance computation.
We visualize them by
embedding the images as points in the plane using MDS.  The results
is shown in \re{Figure~\ref{fig:ells_dists} (top left)}.  
Here the GW distance behaves as expected meaning that the
isometrical pairs are found and located close to each other---the
small visible distances between them result from the chosen
discretizations.  
Up to this expected doubling, we essentially obtain
a triangle, whose corners correspond to the smallest as well as the
largest (isotopic) elliptical disk and the most anisotropic elliptical
disk.  The 4950 pairwise GW distances of this toy
example have been computed in 159.38~minutes.

\re{As mentioned in Remark~\ref{rem:ref_space},
the quality of the approximation of GW by LGW and thus gLGW strongly depends
on the chosen reference space $\SS$.
The choice of $\SS$ is especially crucial
since the minimization over $\Pio(\SS,\XX)$ and $\Pio(\SS, \YY)$ in \eqref{eq:gLGW}
is numerical intractable, 
and fixed optimal plans $\pi_\SS^\XX$ and $\pi_\SS^\YY$ are used instead.
In the sense of Lemma~\ref{lem:LGW_estimates}, 
natural choices for $\SS$ are circular or elliptical disks,
but we also study uniform distributions on squares, triangles, lines 
as well as  composed and non-uniform references.
The employed references are shown in Table~\ref{tab:1}.

To visually compare the approximation quality of $\gLGW_{\SS_i}$
for the considered references $\SS_i$, 
the computed distances are again embedded using MDS,
see Figure~\ref{fig:ells_dists}.
For more quantitative comparisons,
we use the Mean Relative Error (MRE)
and the Pearson Correlation Coefficient (PCC), 
which has been suggested in \cite{VVFCC21} to compare distances.
The computed values are recorded in Table~\ref{tab:1}.
The impact of the different references is clearly visible
raising again the question about a good reference.  

The poorest performances correspond to the circular disk $\SS_1$ and the square $\SS_2$, which have a more regular shape than the others. 
Notice that, 
for any measure-preserving isometry $I:S_i \to S_i$ 
and any optimal plan $\pi_{\SS_i}^{\XX_k} \in \Pio(\SS_i,\XX_k)$,
we have 
\begin{equation}
    \label{eq:bary-plan}
    (I,\id)_\# \pi_{\SS_i}^{\XX_k} \in \Pio(\SS_i,\XX_k);
\end{equation}
so the larger the number of measure-preserving isometries of $\SS_i$, 
the larger the cardinality of $\Pio(\SS_i,\XX_k)$.
Since the square $\SS_2$ is invariant under 3 rotations and 4 reflections
and the circular disk $\SS_1$ nearly under arbitrary rotations and reflections,
the minimization over all optimal plans in \eqref{eq:gLGW}
cannot be neglected any more. 
This issue can also be observed numerically
by examining the computed plans
between the references $\SS_1$, $\SS_2$ and the given mm-spaces $\XX_k$
in Figure~\ref{fig:ref_and_plans}.
Considering the first two columns, 
we notice that
the mass that is transported to the semi-minor axes of the first target
is transported to the semi-major axes of the second target.
Heuristically. $\gLGW_{\SS_i}$ is  
the evaluation of the GW objective in \eqref{eq:GW}
with respect to the plan
\begin{equation*}
\pi \coloneqq \Bigl(\mathcal{T}_{\pi_{\SS_i}^{\XX_k}},\mathcal{T}_{\pi_{\SS_i}^{\XX_\ell}}\Bigr)_\# \sigma_k,
\end{equation*}
where $\mathcal T_{\bullet}$ 
is the generalized barycentric projection in \eqref{eq:gen-bary-proj}.
(Notice that $\pi$ may not satisfy the marginal constraints.)
In this specific instance,
the resulting plan $\pi$ essentially couples the semi-minor axis of the first target with the semi-major axis of the second target,
which is clearly not optimal in the GW sense. 
The consequence of this matching issue is that
$\gLGW_{\SS_1}$ and $\gLGW_{\SS_2}$ cannot recognize
the isometric pairs, which is well reflected by the MDS embedding in Figure~\ref{fig:ells_dists}. 
Although the results of $\SS_3$ are slightly better,
the non-uniform distribution on the square is not able to resolve
this issue numerically.
For the remaining references,
this problem does not occur
since there exist no isomorphic self-couplings or
the self-couplings corresponds to the self-couplings
of the target $\XX_k$---%
rotation by 180° and reflections along the semi-major and semi-minor axes.

\begin{figure}[t]
  \centering
  \includegraphics[width=\linewidth]{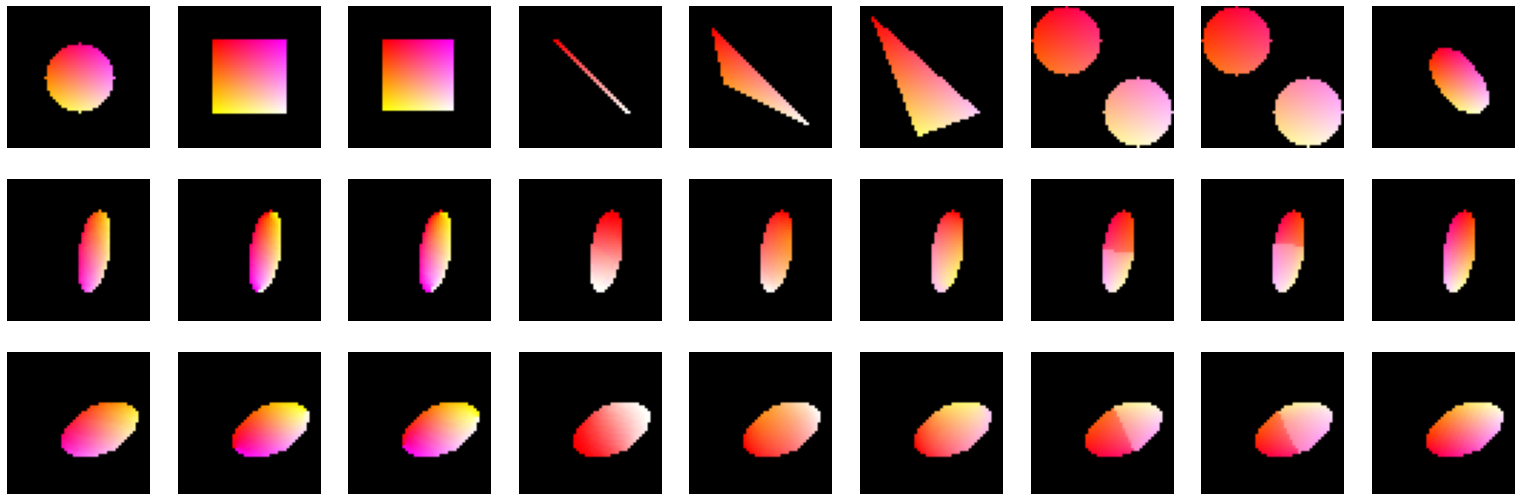}
  \caption{\re{Examples of optimal GW plans between
    the reference mm-spaces (top) and two elliptical disks (middle, bottom) from
    Figure~\ref{fig:ells_data}. The color indicates the mass
    transport from $\SS_i$ to $\XX_k$ that is used for
    the barycentric projection.}}
  \label{fig:ref_and_plans}
\end{figure}

The computation of the 4950 pairwise gLGW distances 
only requires 100 GW transport plans; 
therefore we obtain significant speed-ups in term of computation time.
Considering the qualitative and quantitative results
in Figure~\ref{fig:ells_dists} and Table~\ref{tab:1},
we notice that 
the specific computation time and the MRE strongly depend on
the number support points in the reference space.
The effect on the computation time is clear
since optimal transport plans between spaces with less points
can be calculated faster.
The effect on the MRE is less obvious
and seems to depend on the approximation
of $\XX_k$ by
\begin{equation*}
    \Bigl(\mathcal T_{\pi_{\SS_i}^{\XX_k}} \Bigr)_\# \sigma_k.
\end{equation*}
Heuristically,
the approximation becomes better and the MRE smaller,
if the number of non-zero points in the reference is increased.
Although the MRE with respect to the measure on the line $\SS_4$
is large, $\gLGW_{\SS_4}$ is well correlated to GW.

On the basis of the numerical experiments,
a good reference measure is characterized by
\begin{itemize}
    \item the number of isomorphic GW self-couplings
    (the less, the better) and
    \item the number of non-zero points
    (comparable to the number in the target spaces).
\end{itemize}
Finally, the elliptical disk reference $\SS_9$,
which is close in the GW distance to all given $\XX_k$,
and which should be a good reference in the sense of Remark~\ref{rem:ref_space},
behaves as expected and give excellent results.
}

\subsection{Gromov--Wasserstein in 2D shape analysis}
\label{sec:grom-wass-2d}

\begin{figure}
  \centering
  \includegraphics[width=\linewidth]{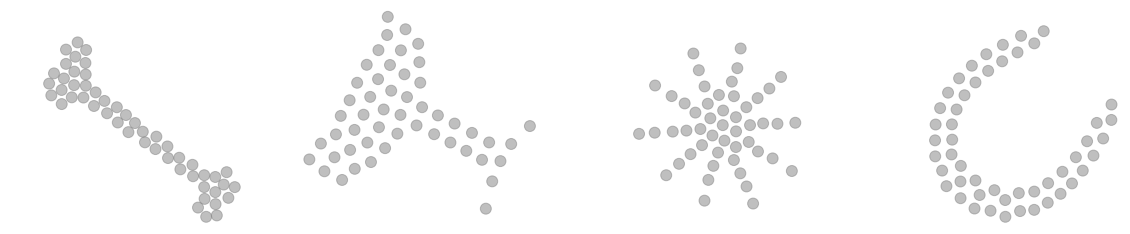}
  \caption{One example of each class (bone, goblet, star, horseshoe)
    of the employed 2D shape dataset.}
  \label{fig:2d-data}
\end{figure}

Next, we apply the GW distance and its linear form to
distinguish different 2D shapes from each other.  In this numerical
experiment, we use the publicly available database \cite{shapes2d}
embracing over 1\,200 shapes in 70 shape classes.  
For our example, we
select 20 shapes of the classes bone, goblet, star, and horseshoe,
respectively, \re{so that we obtain $80$ shapes in total}.  
The shapes are stored as black and white images of
different sizes, where the white pixels corresponds to the objects.
To speed-up the computations, each images is approximated by a point
measure $\mu_k$ consisting of 50 points and uniform weights.  For this
preprocessing step, we use the dithering technique in
\cite{TSGSW11}, see also \cite{EGNS2021}.
\re{All measures are} randomly rotated yielding 80
mm-spaces $\XX_k = ([-1,1]^2, \|\cdot\|, \mu_k) $.  A
preprocessed example of each class is shown in
Figure~\ref{fig:2d-data}.

\begin{figure}
  \centering
  \includegraphics[width=0.25\linewidth]{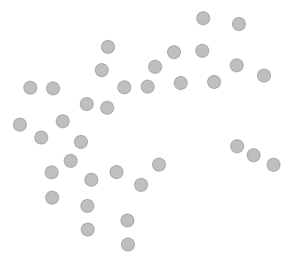}
  \caption{Embedding of the computed barycenter $\SS$ into $\R^2$
    using MDS.}
  \label{fig:2d-bary}
\end{figure}

The performance of the generalized linear GW distance \re{again}
depends on the selection of an appropriate reference space
$\SS$.  
\re{
As discussed in Remark \ref{rem:ref_space},
a barycenter of $\XX_1, \dots, \XX_{80}$ would be a natural choice.
However, 
minimizing \eqref{eq:bary-num} with respect to $80$ inputs is numerically challenging.
Considering the current state-of-the-art algorithm in \cite{PCS2016},
which is based on a blockwise coordinate descent,
we have to compute an optimal GW plan for every input $\XX_k$ per iteration;
so the barycenter computation completely counteracts the computational speed-up by gLGW.
To overcome this issue, we may either approximate the barycenter
by performing only a few iterations 
or exploit that the mm-spaces within the different classes are already 
close to each other.
Following the second approach,
we choose a representative for each of the four classes 
to ensure that all main features are covered
and compute a barycenter with 35 points and
uniform weights.
The employed barycenter is shown in Figure~\ref{fig:2d-bary}.  
Using the POT package,
the computation takes 1.60 seconds.
Since the reference $\SS$ has less points than the spaces $\XX_k$, 
the barycentric projection \eqref{eq:baryproj} is indeed a mapping to
weighted means.}

\begin{figure}
  \centering
    \includegraphics{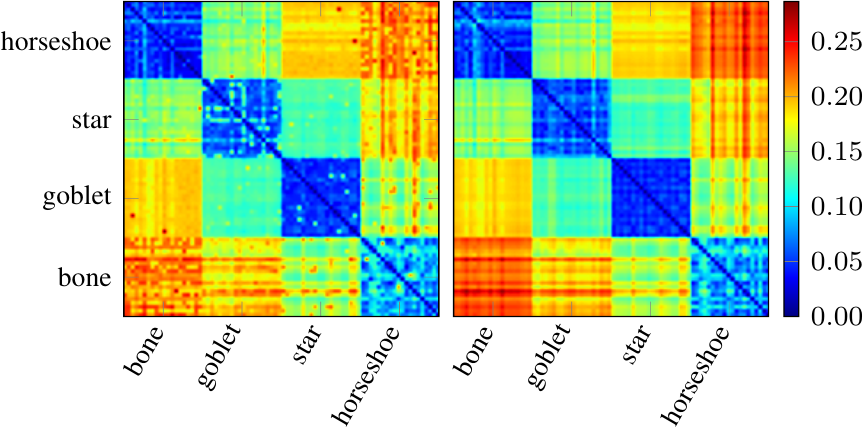}
  \caption{Pairwise GW distances (left) and 
    generalized linear GW distances for a barycenter
    reference $\SS$.  All images share the same color coding.}
  \label{fig:2d-gw-dist}
\end{figure}

On the basis of the chosen reference $\SS$, we now compute the
pairwise GW distances (11.34 seconds) 
and  
the pairwise gLGW distances (0.72 seconds). 
\re{Even with barycenter computation,
gLGW gives a significant speed-up against GW.}
The results are shown in
Figure~\ref{fig:2d-gw-dist}.  Notice that the shapes in the horseshoe
class significantly differ between each other explaining the greater
distances than in other classes.  \re{In this example, the computed GW and gLGW distances are visually comparable.  A more quantitative comparison is given below.}

\begin{figure}
  \centering
  \includegraphics{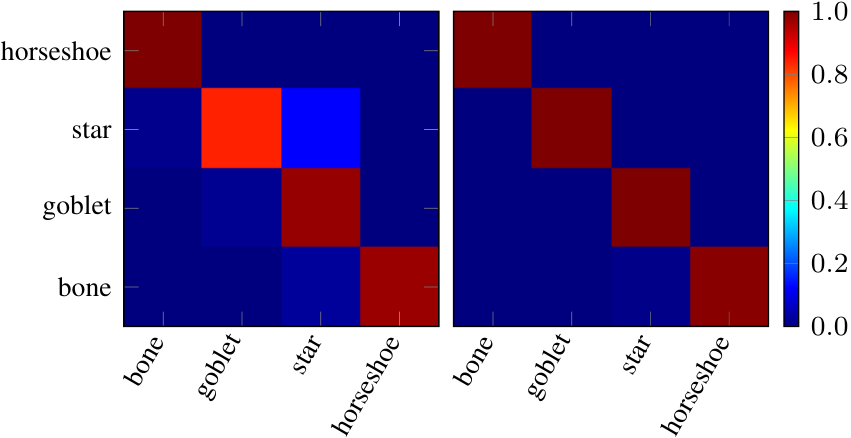}  
  \caption{Confusion matrix for the GW distances
    (left) and the generalized linear GW distances
    with barycenter reference (right).}
  \label{fig:2d-conf-mat}
\end{figure}

Considering the results, it seems reasonable to use a nearest neighbor
classification to distinguish the different classes with
respect to some representatives.  Numerically, this concept may be
verified by computing a \emph{confusion matrix} consisting of the
probabilities to classify an instance of a class as another class.
For this, we rely on \cite[§~8.2]{memoli2011gromov}, where the
confusion matrix is estimated by randomly choosing a representative
for each class and then classifying all other shapes $\XX_k$,
$k = 1,\dots,80$, with respect to the nearest representative.  This
classification task is then repeated 10\,000 times.  The confusion
matrix of $\GW$ and $\gLGW$ are shown in
Figure~\ref{fig:2d-conf-mat}.  \re{Interestingly, gLGW performs slightly
better than GW.}

\begin{table}
    \caption{\re{Quantitative comparisons of the 10-fold cross-validation. The recorded values correspond to the mean over all 10 training and classification tasks. 
    The SVM has been trained based on gLGW, where barycenters of different sizes have been used.}}
    \label{tab:cv-2d}
    \centering
    \re{
    \begin{tabular}{lccccc}
        \toprule
        points in barycenter & 15 & 25 & 35 & 45 & 55 
        \\ \midrule
        mean accuracy         
        & 0.9750 & 0.9625 & 0.9875 & 0.975 & 0.9875
        \\
        mean MRE
        & 0.2449 & 0.3040 & 0.1623 & 0.1826 & 0.1709
        \\
        mean PCC
        & 0.8042 & 0.7875 & 0.8936 & 0.8836 & 0.8904
        \\ \bottomrule
    \end{tabular}}
\end{table}

\re{
The nearest neighbor classification already shows
that the distinctiveness of gLGW is comparable to GW.
To provide a more quantitative study,
we combine gLGW with a support vector machine (SVM),
see for instance \cite{Ste15} and references therein.
Following the approach in \cite{TCTLF19},
we employ the kernels $\exp(-\alpha\GW)$ and $\exp(-\alpha\gLGW_\SS)$
with $\alpha > 0$
although these might not be positive definite.
We obtain the best performance for $\alpha \coloneqq 10$.
Moreover, we apply a 10-fold cross-validation.
For this, we divide the given dataset $\XX_1, \dots, \XX_{80}$
with respect to the classes into 10 disjoint subsets.
In each iteration of the cross-validation,
we train the SVM based on 9/10 of the data and use the 
remaining 1/10 data as test set.
Further,
the employed barycenter $\SS$ is computed anew from random representatives of the four classes with respect to the current training set.
After each training,
we compute the empirical success rate (accuracy) of the classification according to the current test set as well as the MRE and PCC between GW and gLGW. 
The means over all 10 cross-validation steps 
for different sizes of the barycenter are recorded in Table~\ref{tab:cv-2d}.
Note that the SVM with respect to the GW distance achieves a perfect accuracy score of one.
Using gLGW, we encounter up to three misclassifications over all
10 cross-validation steps in total.
Considering the mean classification accuracy, we may reduce the size of the barycenter to 15 points, which additionally speeds up the barycenter and gLGW computations. 
The MRE and PCC are improved for higher numbers.  
Both reach their optimum at around 35 points, which corresponds to the former given qualitative results. 
The appropriate number of points in the barycenter thus mainly depends on the application.  If we are interested in classification, we may choose less points;  if we are interested in pairwise GW approximations, we require more points. 
}

\subsection{Gromov--Wasserstein in 3D shape analysis}

The GW distance traces back to the comparison and
matching of 3D shapes, which we take up in our final numerical
example.  Analytically, a 3D shape is a two-dimensional submanifold of
$\R^3$ that may have a boundary.  3D shapes can be interpreted as
mm-spaces $\XX = (X, d_X, \mu)$, where $X$ is a surface of the shape, where
$d_X$ corresponds to the length of the geodesics between two points,
and where $\mu$ is some measure.

In practice, 3D shapes are usually triangulated and thus realized by a
net of triangles.  To handle them numerically, we approximate them by
a discrete mm-space $\XX = (X, d_X, \mu)$.  The vertices of the net
become the discrete points in $X$.  To approximate the geodesic
distance on $X$, a \emph{weighted graph} $G = (X, E)$ consisting of
all vertices $X$ and all edges $E$ of the triangulation may be used,
where the edges are weighted by the Euclidean distance between the
corresponding vertices.  The geodesic distance between two vertices
may now be approximated by the length of the shortest path between
these vertices. This distance can be computed by the Dijkstra algorithm from
the NetworkX package \cite{networkx}. \re{The probability measure $\mu$ may
be used to incorporate additional information of the shapes.}

\begin{figure}
    \centering
    \includegraphics[width = 0.16\linewidth]{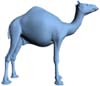}\quad
    \includegraphics[width = 0.16\linewidth]{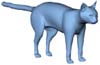}\quad
    \includegraphics[width = 0.16\linewidth]{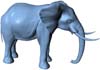}\quad
    \includegraphics[width = 0.12\linewidth]{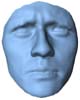}\\
    \includegraphics[width = 0.07\linewidth]{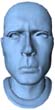}\quad
    \includegraphics[width = 0.16\linewidth]{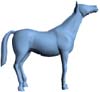}\quad
    \includegraphics[width = 0.16\linewidth]{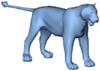}\quad
    \caption{One example of each class (camel, cat, elephant, face,
      head, horse, lion) of the employed 3D shape dataset.}
    \label{fig:3d_meshes}
\end{figure}

In this example, we consider 3D shapes of the publicly available
database \cite{mesh3d_animals}, which has already been used by Mémoli
\cite[§~8.2]{memoli2011gromov} in the context of
GW distances. We use a similar setting to make the results comparable.  
As dataset for the experiment, we
choose 3D shapes corresponding to the animals camel, cat, elephant,
horse, and lion as well as to a human face and head.  Each object is
shown in 10 to 11 different poses totaling to 73 shapes.
Figure~\ref{fig:3d_meshes} shows one example pose of each object.
Every object is provided by a triangulation consisting of up to 43.000
vertices and up to 130.000 triangles.

Since the discrete mm-spaces of the full triangulations consist of
too many points for our purpose, the 3D shapes are preprocessed by a
two-step approximation similar to \cite{memoli2011gromov}.
\begin{enumerate}
\item Starting from a given triangulation with vertices $X$, 
we first reduce $X$ to a set $\widetilde X$ consisting of 4\,000
vertices.  The first vertex is hereby chosen randomly and is
sequentially followed by the points with the largest Euclidean distance 
to the already chosen points.  This selection rule is also known as
the furthest point procedure.  
\item The set $\widetilde X$ is reduced
further and an appropriate measure $\mu$ is constructed.  
For this, we again apply the furthest point procedure to reduce $\widetilde X$ to a
subset $\widehat X$ consisting of 50 points, but this time with respect to the discrete
geodesic distance $d_X$ calculated using a weighted graph as explained
above. Then we endow $\widehat X$
with a discrete probability measure, where the mass at $x \in \widehat
X$ is proportional to the amount of closest neighbors within $\widetilde
X$ with respect to the original geodesic distance $d_X$.
In other words, we compute the Voronoi diagram of $\widetilde X$ to the
points $\widehat X$ with respect to $d_X$ and count the members of
each Voronoi cell.  Repeating this procedure for every given 3D shape,
we end up with 73 discrete mm-spaces $\XX_k = (\widehat X_k, d_{X_k},
\mu_k)$, $k = 1, \dots, 73$.
\end{enumerate}

Since the distance $d_X$ of the constructed mm-spaces $\XX_k$ are discrete
geodesic distances, which are restricted to the points in
$\widehat X_k$, the barycentric projection
\eqref{eq:gen-bary-proj} has the form
\begin{equation*}
  \mathcal T_{\pi_\SS^{\XX_k}}(s)
  \coloneqq \argmin_{x_0 \in \widehat X_k} \sum_{x \in \widehat X_k}
  \pi_\SS^{\XX_k} (\re{\{(s,x)\}}) \, {d_{X_k}^2(x_0,x)},
\end{equation*}
where $\pi_\SS^{\XX_k}$ is the chosen optimal GW
transport plan.

\begin{figure*}
    \centering
    \includegraphics{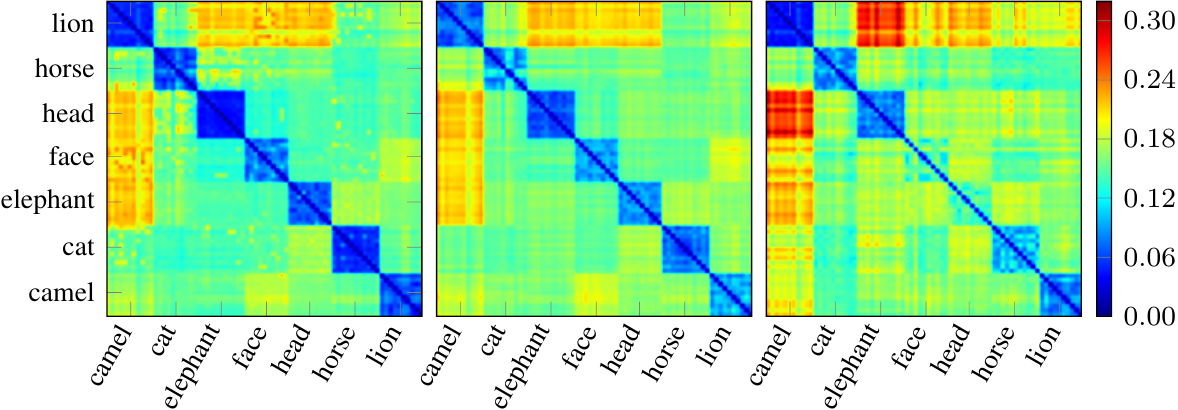}
    \caption{Pairwise GW distances (left), 
      generalized linear GW distances for a
      barycenter reference $\SS$ (middle), and generalized
      linear GW distances with reference
      $\SS = \XX_{10}$ (camel) (right).  All images share the same
      color coding.}
    \label{fig:pw_gw_vs_lgw_memoli}
\end{figure*}

The pairwise GW distances (14.29 seconds) and the generalized
linear GW distances (0.82 seconds) for the 3D shape
dataset with respect to two different reference measures 
are shown in Figure~\ref{fig:pw_gw_vs_lgw_memoli}.  
\re{Similarly to the 2D shape example in Section~\ref{sec:grom-wass-2d},
one of the considered references is a GW barycenter. 
To speed-up the computations, we again choose one representative of each class.
To compute the barycenter
with 50 points corresponding to a uniform distribution,
the POT package needs around 4.67 seconds.
Considering the middle image of Figure~\ref{fig:pw_gw_vs_lgw_memoli},
we notice that gLGW is again comparable to GW,
i.e.\ the different classes are clearly identifiable.  
Analogously to the previous numerical examples,
and as indicated by Remark~\ref{rem:ref_space},
the barycenter gives excellent results.
Since the computation of the barycenter is, however, numerically costly,
we secondly use the given reduced 3D shape $\SS \coloneqq \XX_{10}$ (camel)
as reference.  
Here, the quality of the GW approximation illustrated in the right-hand side image of Figure~\ref{fig:pw_gw_vs_lgw_memoli} is more diverse.
On the one side,
the approximation of the GW distances inside the camel class 
is nearly perfect;
on the other side, the approximation outside the camel class loses in quality.
Especially, the head and face classes are affected.}

\begin{figure*}
    \centering
    \includegraphics{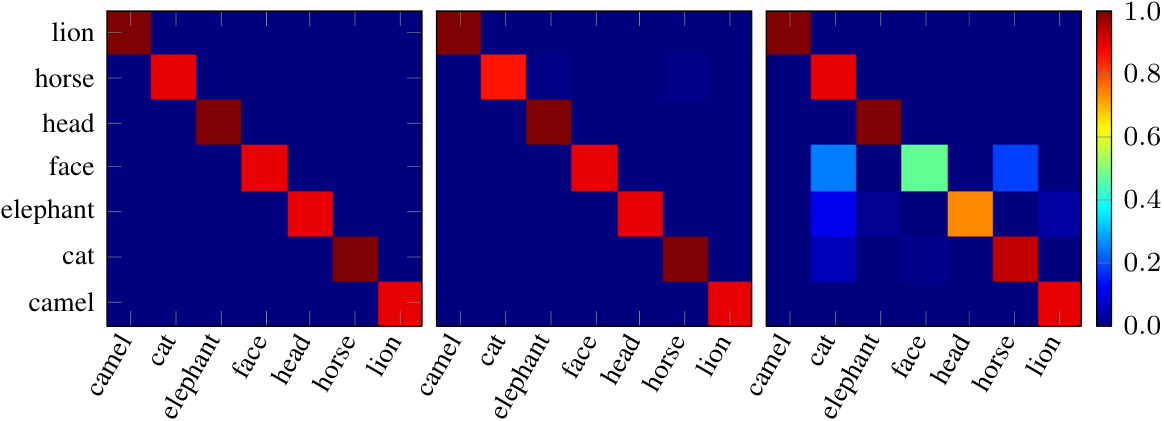}
    \caption{Confusion matrices for the GW distances
      (left), the generalized linear GW distance with
      barycenter reference (middle), and the generalized linear
      GW distance with reference $\SS =
      \XX_{10}$ (right).}
    \label{fig:conf_mat}
\end{figure*}

To evaluate the classification quality of the $\GW$ and $\gLGW$
distances, we compute the confusion matrix $C$ for both
distances as in Section~\ref{sec:grom-wass-2d}.  The confusion matrix
$C$ consists of the probabilities to classify a 3D shape within one
class (camel, cat, elephant, face, head, horse, lion) to another
class.  Following again \cite[§~8.2]{memoli2011gromov}, for this
purpose, we first randomly chose a representative for each class and
then classify all 3D shapes $\XX_k$, $k = 1,\dots,73$, with respect to
the nearest representative.  This classification task is repeated
10\,000 times.  The result is shown in Figure~\ref{fig:conf_mat}.
\re{Considering the first two confusion matrices, 
we notice that the results for GW (left) and gLGW
with barycenter reference (middle) nearly coincide.
The classification based on gLGW with reference $\XX_{10}$ (right) performs slightly worse.
Camel $\XX_{10}$ allow the classification of all four-legged animals.
The serious misclassifications occur for the head and, especially, the face class, which is not astonishing since the geometry of $\XX_{10}$
is quite differently from these two classes.
 }

\begin{table}
    \caption{\re{Mean MRE achieved during a 10-fold cross-validation.  For each data splitting a new barycenter is computed.  The experiment is repeated for different numbers of support points in the reduced graph and in the barycenter.  The last row records the mean MRE for the camel reference $\XX_{10}$.}}
    \label{tab:mre-3d}
    \centering
    \re{
    \begin{tabular}{ccccc}
        \toprule
        & \multicolumn{4}{c}{points in graph}
        \\ \cmidrule{2-5}
        points in barycenter & 25 & 50 & 75 & 100 
        \\ \midrule
        25 & 0.2339 & 0.1659 & 0.1687 & 0.1745 \\
        50 & 0.1975 & 0.1598 & 0.1486 & 0.1425 \\
        75 & 0.2408 & 0.1884 & 0.1798 & 0.1586 \\
        100 & 0.2417 & 0.2083 & 0.1799 & 0.1641 \\
        \midrule
        camel $\XX_{10}$ & 0.2204 & 0.2124 & 0.2179 & 0.2618
        \\ \bottomrule
    \end{tabular}}
\end{table}

\begin{table}
    \caption{\re{Mean PCC achieved during a 10-fold cross-validation.  For each data splitting a new barycenter is computed.  The experiment is repeated for different numbers of support points in the reduced graph and in the barycenter.  The last row records the mean PCC for the camel reference $\XX_{10}$.}}
    \label{tab:pcc-3d}
    \centering
    \re{
    \begin{tabular}{ccccc}
        \toprule
        & \multicolumn{4}{c}{points in graph}
        \\ \cmidrule{2-5}
        points in barycenter & 25 & 50 & 75 & 100 
        \\ \midrule
        25 & 0.7577 & 0.8225 & 0.8452 & 0.8282 \\
        50 & 0.8686 & 0.9251 & 0.9268 & 0.9318 \\
        75 & 0.8537 & 0.9055 & 0.9156 & 0.9302 \\
        100 & 0.8769 & 0.9080 & 0.9212 & 0.9370 \\
        \midrule
        camel $\XX_{10}$ & 0.7300 & 0.8428 & 0.8511 & 0.8496 
        \\ \bottomrule
    \end{tabular}}
\end{table}

\begin{table}
    \caption{\re{Mean accuracy achieved during a 10-fold cross-validation.
    For each data splitting a new barycenter is computed.  The experiment is repeated for different numbers of support points in the reduced graph and in the barycenter.  The last row records the mean accuracy for the camel reference $\XX_{10}$.}}
    \label{tab:acc-3d}
    \centering
    \re{
    \begin{tabular}{ccccc}
        \toprule
        & \multicolumn{4}{c}{points in graph}\\
        \cmidrule{2-5}
        points in barycenter & 25 & 50 & 75 & 100 \\
        \midrule
        25 & 0.9482 & 1.0000 & 0.9875 & 0.9714 \\
        50 & 0.9875 & 0.9857 & 1.0000 & 0.9857 \\
        75 & 0.9607 & 1.0000 & 1.0000 & 0.9857 \\
        100 & 1.0000 & 1.0000 & 0.9857 & 0.9750 \\
        \midrule
        camel $\XX_{10}$ & 0.8535 & 0.9250 & 0.8964 & 0.9107
        \\ \bottomrule
    \end{tabular}}
\end{table}

\re{
Similarly to the previous experiment, 
we train a SVM with respect to 
$\exp(-\alpha\GW)$ and $\exp(-\alpha\gLGW_{\SS})$
with $\alpha \coloneqq 10$. 
The parameter choice $\alpha = 10$ has again performed best.
During the applied 10-fold cross-validation,
a new barycenter is computed for every data splitting 
from one random representative of each class in the training set.
Moreover, we repeat the cross-validation for different sizes of the reduced graphs, i.e. numbers of support points in $\widehat X_k$,
and different sizes of the barycenter.
The resulting MRE and PCC are recorded in Table~\ref{tab:mre-3d} and \ref{tab:pcc-3d}.
The best performance is obtained for the constellations of 50 points 
in the barycenter and 50 or more points in the reduced graph.
Comparing both tables,
we notice that the gLGW procedure performs best 
if the number of support points in the barycenter
is less or equal the number of
support points in the target spaces. 
As comparison, the last row in each table records the
performance of gLGW with $\SS = \XX_{10}$,
where the number of points in the reference and the reduced graph coincide,
and where the datum $\XX_{10}$ has been removed from the training and testing subsets.
Moreover,
the SVM based on GW
archives perfect accuracy scores of one. 
The accuracy with respect to gLGW
is recorded in Table~\ref{tab:acc-3d}. 
Although there are some misclassifications with respect
to the head and face classes, 
the accuracy score for the camel reference $\SS=\XX_{10}$ is high.
The numerical experiments show that the classification by the
SVM with gLGW is powerful
even for non-optimal references.}

\section{Conclusions}\label{sec:conclusions}
We proposed a linear version of the GW distance
that was inspired by a generalized version of the linear Wasserstein distance.
As the latter one, the approach appears to be efficient in applications,
where pairwise distances of a larger amount of measures are of interest.
We gave three examples indicating that 
our linear version of the GW distance gives reasonable approximations
and circumvents the heavy computation of all pairwise distances.
In contrast to Wasserstein distances, 
the mathematics behind GW distances 
is not well-examined so far and there are plenty of open problems
which could be tackled in the future.
For example, in generalized version of LOT, it would also be possible 
to use the concept of weak optimal transport \cite{GRST2017}. 
This approach was neither
considered for gLOT nor for gLGW so far. Further, multimarginals may be addressed, see \cite{BLNS2021}.
Finally, we are interested in further \re{applications
in the context of shape and graph analysis.  
More precisely, we would like to incorporate our generalized linear Gromov--Wasserstein distance into existing shape and graph classification approaches exploiting feature spaces, annotations, and deep learning.}

\section*{Acknowledgment}
The funding by the German Research Foundation (DFG) within the RTG
2433 DAEDALUS and by the BMBF project ``VI-Screen'' (13N15754) is
gratefully acknowledged.  Further, the authors would like to thank
Johannes von Lindheim for valuable discussions and for assisting with
numerical implementations \re{as well as the anonymous reviewers for their valuable comments and suggestions to improve the manuscript
and to strengthen the numerical simulations.}

{\appendices
\section{Proofs}

\subsection{Proof of Proposition \ref{prop_1}}
We have that $\tilde \mu \in \mathcal P_2(\R^d)$ since by Jensen's inequality
\begin{align*}
  &\int_{\R^d} \| x \|^2 \dx \tilde \mu(x)
  \\
  &= \int_{\R^d} \| \mathcal T_{\pi_\sigma^\mu} (s) \|^2 \dx \sigma(s)
    = \int_{\R^d} \biggl\| \int_{\R^d} x \dx \mathcal \pi_s(x) \biggr\|^2 \dx \sigma(s)
  \\
  &\le \int_{\R^d} \int_{\R^d}
    \| x \|^2 \dx \pi_s(x) \dx \sigma(s)
    = 
    \int_{\R^d \times \R^d} \| x \|^2 \dx \pi_\sigma^\mu
    (s,x)
  \\
  &= \int_{\R^d} \|x\|^2 \dx \mu(x) < \infty.
\end{align*}
Let $\pi_\sigma^{\tilde \mu}$ be an optimal transport plan 
with respect to $W(\sigma,\tilde \mu)$.
By the dual formulation of the optimal transport problem, 
see \cite[Thm~4.2]{ABS21} and \cite[Thm~5.10]{villani2008optimal},
we know that
\begin{align}
&\int_{\R^d \times \R^d} \frac12 \| s - x \|^2 \dx \pi_\sigma^{\tilde \mu} (s,x)
  \\
  &= 
\sup_{\phi \in L^1_\sigma(\R^d)} \biggl\{ \int_{\R^d} \phi(s) \dx \sigma(s) + \int_{\R^d} \phi^c(x) \dx \tilde \mu(x) \biggr\},
\end{align}
where $\phi^c$ denotes the $c$-concave function  given by 
$$\phi^c (x) = \inf_{y \in X} \{ \tfrac 12\|x-y\|^2 - \phi(y) \}.$$
To yield a contradiction, assume that $\mathcal T_{\pi_\sigma^\mu}$ is not an optimal transport map. Then 
$$ \tilde \pi \coloneqq (\id, \mathcal  T_{\pi_\sigma^\mu})_\# \sigma$$
is not an optimal transport plan with respect to $W(\sigma,\tilde \mu)$
and
\begin{align}
  &\int_{\R^d} \frac12 \| s - \mathcal T_{\pi_\sigma^\mu}(s) \|^2 \dx \sigma(s) 
\\
  &= 
    \int_{\R^d \times \R^d} \frac12 \| s - x \|^2 \dx \tilde \pi (s,x)  \nonumber \\
  &
    > 
    \sup_{\phi \in L^1_\sigma(\R^d)} \biggl\{ \int_{\R^d} \phi(s) \dx \sigma(s) + \int_{\R^d} \phi^c(x) \dx \tilde \mu(x) \biggr\}. 
    \label{ass}
\end{align}
Now we obtain for the optimal transport plan $\pi_\sigma^\mu$ of $W(\sigma,\mu)$ that
\begin{align*}
  &\int_{\R^d \times \R^d} \| s - x \|^2 \dx \pi_\sigma^\mu (s,x)
  \\
  &= \int_{\R^d} \| s \|^2 \dx \sigma(s)
    -2 \int_\R^d \int_\R^d \langle s,x \rangle \dx
    \pi_s(x) \dx \sigma(s)
  \\
  &\qquad
    + \int_{\R^d} \| x \|^2 \dx \mu(x)
  \\
  &= \int_{\R^d} \| s \|^2 \dx \sigma(s)
    -2 \int_\R^d \langle s, \mathcal T_{\pi_\sigma^\mu}(s)
    \rangle \dx \sigma(s)
  \\
  &\qquad
    + \int_{\R^d} \| x \|^2 \dx \mu(x)
  \\
  &= \int_{\R^d} \| s - \mathcal T_{\pi_\sigma^\mu}(s) \|^2 \dx \sigma(s)
    + \int_{\R^d} \| x \|^2 \dx \mu(x)
  \\
  &\qquad
    - \int_{\R^d} \| x \|^2 \dx \tilde \mu(x)
\end{align*}
and by \eqref{ass} further
\begin{align}
  &\int_{\R^d \times \R^d} \frac 12 \| s - x \|^2
    \dx \pi_\sigma^\mu (s,x)
  \\
  &> \int_{\R^d} \frac12 \| x \|^2 \dx \mu(x)
    - \int_{\R^d} \frac12 \| x \|^2 \dx \tilde \mu(x)
  \\
  &\qquad
    + \sup_{\phi \in L^1_\sigma(\R^d)} 
    \biggl\{ \int_{\R^d} \phi(s) \dx
    \sigma(s) + \int_{\R^d} \phi^c(x) \dx \tilde
    \mu(x) \biggr\}
  \\
  &= \sup_{\phi \in L^1_\sigma(\R^d)} 
    \biggl\{\int_{\R^d} \frac12 \| x \|^2 \dx \mu(x)
    - \int_{\R^d} \frac12 \| x \|^2 \dx \tilde \mu(x)
  \\
  &\qquad+  \int_{\R^d} \phi(s) \dx
    \sigma(s) + \int_{\R^d} \phi^c(x) \dx \tilde
    \mu(x) \biggr\}
  \\
  &= \sup_{\phi \in L^1_\sigma(\R^d)} 
    \biggl\{\int_{\R^d} \frac12 \| x \|^2 \dx \mu(x)
    +  \int_{\R^d} \phi(s) \dx
    \sigma(s)
  \\
  &\qquad+ \int_{\R^d} \phi^c(x) - \frac12 \| x \|^2 \dx \tilde
    \mu(x) \biggr\},
    \label{eq:prim-dual}
\end{align}
where $L^1_\sigma(\R^d)$ is the space of functions which absolute values are integrable with respect to $\sigma$.
Since $\phi^c$ is $c$-concave, we know that $h \coloneqq \phi^c - \tfrac12 \|\cdot\|^2$ is concave, 
see \cite[Lect~4.4]{ABS21}.
Thus, Jensen's inequality implies
\begin{align}
  &\int_{\R^d} \phi^c(x) - \frac12 \| x \|^2 \dx \tilde \mu(x)
\\
  &= 
    \int_{\R^d} \phi^c(\mathcal T_{\pi_\sigma^\mu}(s)) - \frac12
    \| \mathcal T_{\pi_\sigma^\mu}(s) \|^2 \dx \sigma(s)
    \\
    &= \int_{\R^d} h\left(\mathcal  T_{\pi_\sigma^\mu}(s) \right) \dx \sigma(x)
  \\
  &= 
    \int_{\R^d} h\left( \int_{\R^d}  x \dx \pi_s(x) \right) \dx \sigma(s)
  \\
  &
    \ge \int_{\R^d} \int_{\R^d} h(x) \dx \pi_s(x)  \dx \sigma(s)
  \\
  &= 
    \int_{\R^d} \int_{\R^d} \phi^c(x) - \frac12 \| x \|^2 \dx 
    \pi_s(x) \dx \sigma(s)
  \\
  &=
    \int_{\R^d} \phi^c(x) \dx \mu(x)
    - \int_{\R^d} \frac12 \| x \|^2 \dx \mu(x).
    \label{eq:c-conc-jens}
\end{align}
Inserting \eqref{eq:c-conc-jens} into \eqref{eq:prim-dual}, we obtain
\begin{align*}
  &\int_{\R^d \times \R^d} \frac12 \| s - x \|^2 \dx \pi_\sigma^\mu (s,x)
  \\
  &> \sup_{\phi \in L^1_\sigma(\R^d)} 
  \biggl\{ \int_{\R^d} \phi(s) \dx \sigma(s) + \int_{\R^d} \phi^c(x) \dx \mu(x) \biggr\},
\end{align*}
which contradicts the optimality of $\pi_\sigma^\mu$. 
\hfill $\Box$

\subsection{Proof of Proposition \ref{thm:exp-gw-geo}}

To compute the distance \eqref{dist_final},
recall that the geodesics related to
$\pi_\SS^\XX \in \Pio(\SS, \XX)$ and $\pi_\SS^\YY \in \Pio(\SS, \YY)$
are mapped to the tangents
\begin{align*}
  g
  \coloneqq F_\SS(\pi_\SS^\XX)
  &= d_X - d_S
  &&(\text{acting on} \;
    \TT_g \coloneqq (S \times X, d_S, \pi_\SS^\XX)),
  \\
  h
  \coloneqq F_\SS(\pi_\SS^\YY)
  &= d_Y - d_S
  &&(\text{acting on} \;
    \TT_h \coloneqq (S \times Y, d_S, \pi_\SS^\YY)),
\end{align*}
where $\SS \sim \TT_g \sim \TT_h$.
We next characterize the plans $\pi \in \Pio(\TT_g, \TT_h)$
occurring in the definition of $\GW_\SS$ in \eqref{eq:gw-tan}.
Since $\TT_g$ and $\TT_h$ are equivalent,
each plan  $\pi \in \Pio(\TT_g, \TT_h)$ satisfies
\begin{align*}
  0
  &= \GW(\TT_g, \TT_h)\\
  &=\int_{(S \times X \times S \times Y)^2}
    \hspace{-45pt}
    | d_S(s_1, s_2) - d_S(s_1', s_2') |^2
  \\[-5pt]
  & \hspace{80pt}
    \times \dx \pi(s_1,x,s_2,y) \dx \pi(s_1',x',s_2',y'),
  \\
  & =\int_{(S  \times S)^2}
    \hspace{-20pt}
    | d_S(s_1, s_2) - d_S(s_1', s_2') |^2
    \dx \gamma(s_1,s_2) \dx \gamma(s_1',s_2'),
\end{align*}
where $\gamma \coloneqq P^{13}_\# \pi$.
Thus,
$\gamma$ is an optimal self-coupling of $\SS$ in the GW sense.
As stated in \cite[Lem~10.4]{memoli2011gromov},
each self-coupling has the form $\gamma = (\id,\psi)_\# \sigma$
for some measure-preserving isometry $\psi:S \to S$.
This, however, implies that
the mapping $P^{124}(s_1,x,s_2,y) = (s_1,x,y)$ is $\pi$-almost everywhere invertible
by $(P^{124})^{-1} (s, x, y) = (s, x, \psi(s), y)$.
More precisely,
$(P^{124})^{-1} \circ P^{124}$ is the identity on $\supp( \pi)$.
Therefore,
we have
\begin{equation*}
  \pi = (P^{124})^{-1}_\# \tilde \pi
  \quad\text{with}\quad
  \tilde \pi = P^{124}_\# \pi.
\end{equation*}
Considering the marginals of $\tilde \pi$,
we conclude that
every 4-plan $\pi \in \Pio(\TT_g,\TT_h)$ can be uniquely identified
by the 3-plan $\tilde{\pi} \in \re{\Gamma}_\SS(\pi_\SS^\XX, \pi_\SS^\YY)$,
and vice versa.
This identification finally allows us
to rewrite the metric \eqref{eq:gw-tan} on the tangent space
using the substitution $\pi = (P^{124})^{-1}_\# \tilde \pi$
to obtain
\begin{align*}
  &\GW_\SS(\pi_\SS^\XX, \pi_\SS^\YY)
  \\
  &= \inf_{\pi \in \Pio(\TT_g, \TT_h)}
    \int_{(S \times X \times S \times Y)^2}
    \hspace{-47pt}
    \bigl\lvert d_X(x,x') - d_S(s_1,s_1')
  \\
  & \hspace{100pt}
    - d_Y(y,y') + d_S(s_2,s_2') \bigr\rvert^2
  \\
  &\hspace{100pt}
    \times \dx \pi(s_1,x, s_2,y) \dx \pi(s_1',x', s_2',y'),
  \\
  &= \inf_{\tilde \pi \in \re{\Gamma}_\SS(\pi_\SS^\XX, \pi_\SS^\YY)}
    \int_{(S \times X \times Y)^2}
    \hspace{-35pt}
    \bigl\lvert d_X(x,x') - d_Y(y,y')
  \\[-7pt]
  & \hspace{115pt}
     +d_S(\psi(s), \psi(s')) - d_S(s,s') \bigr\rvert^2
  \\
  &\hspace{115pt}
    \times \dx \tilde \pi(s,x,y) \dx \tilde \pi(s',x',y'),
  \\
  &= \inf_{\tilde \pi \in \re{\Gamma}_\SS(\pi_\SS^\XX, \pi_\SS^\YY)}
    \int_{(S \times X \times Y)^2}
    \hspace{-35pt}
    \lvert d_X(x,x') - d_Y(y,y') \rvert^2
  \\[-7pt]
  &\hspace{132pt}
    \times \dx \tilde \pi(s,x,y) \dx \tilde \pi(s',x',y'),
\end{align*}
which establishes the assertion.
\hfill $\Box$

\re{
\subsection{Proof of Lemma~\ref{lem:LGW_estimates}}
For every $\pi_\SS^\XX \in \Pio(\SS,\XX)$ and $\pi_\SS^\YY \in \Pio(\SS,\YY)$ 
in \eqref{eq:GW_sigma}, 
a three-plan $\pi \in \re{\Gamma}_\SS(\pi_\SS^\XX,\pi_\SS^\YY)$ in \eqref{eq:GW_S} satisfies $P^{23}_\# \pi \in \Pi(\mu,\nu)$. 
Hence $\GW(\XX,\YY) \leq \LGW_\SS(\XX,\YY)$. 
For the upper bound, we consider fixed $\pi_\SS^\XX \in \Pio(\SS,\XX)$ and $\pi_\SS^\YY \in \Pio(\SS,\YY)$ 
in the definition of $\LGW_\SS$ in \eqref{eq:GW_sigma}.
Exploiting the definition of $\GW_\SS$ in \eqref{dist_final} 
and the metric \eqref{eq:gw-tan},
we have
\begin{align*}
\GW_\SS(\pi_{\SS}^\XX,\pi_{\SS}^\YY) &= \GW_\SS(d_X - d_S,d_Y - d_S) \\
&\leq \GW_\SS(d_X - d_S,0) + \GW_\SS(d_Y - d_S,0),
\end{align*}
where $0$ denotes the zero function on $S \times X$.
Using the representatives $\TT_{d_X-d_S} \coloneqq (S \times X, d_X - d_S, \pi_\SS^\XX)$
and $\TT_0 \coloneqq (S \times X, 0, \pi_\SS^\XX)$, 
we further obtain
\begin{align*}
    &\GW_\SS(d_X - d_S,0) \\
    &= \inf \bigl\{\|d_X - d_S\|_{L^2((S \times X)^4, \pi \otimes \pi)} : \pi \in \Pio(\TT_{d_X-d_S},\TT_0) \bigr\}\\
    &=
    \inf_{\pi \in \Pio(\TT_{d_S-d_X}, \TT_0)}
    \int_{(S \times X \times S \times X)^2}
    \hspace{-47pt}
    |d_X(x_1,x_1') - d_S(s_1, s_1')|
    \\
    &\hspace{87pt}
    \times \dx \pi(s_1,x_1,s_2,x_2)
    \dx \pi(s_1',x_1',s_2',x_2')
    \\
    &= \int_{(S \times X)^2}
    \hspace{-23pt}
    |d_X(s_1,x_1) - d_S(s_1', x_1')|
    \dx \pi_\SS^\XX(x_1,s_1)
    \dx \pi_\SS^\XX(x_1',s_1')
    \\
    &= \GW(\SS,\XX)
\end{align*}
since $P^{12}_\# \pi = \pi_\SS^\XX$ for all $\pi \in \Pio(\TT_{d_X-d_S},\TT_0)$.
A similar computation shows $\GW_\SS(d_Y - d_S,0) = \GW(\SS,\YY)$.
Thus  $\LGW_\SS(\XX,\YY) \leq \GW(\SS,\XX) + \GW(\SS,\YY)$ as desired. 
\hfill $\Box$
}}

\bibliographystyle{IEEEtran}
\bibliography{IEEEabrv,references}

\begin{IEEEbiography}[{\includegraphics[width=1in,height=1.25in,clip,keepaspectratio, trim=0pt 35pt 0pt 0pt]{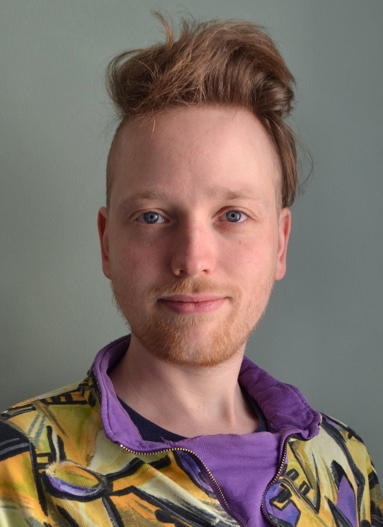}}]{Florian Beier}
studied mathematics at the Technical University Berlin. After recieving his MSc in 2021, he started his PhD under the supervision of Gabriele Steidl. His main research interest is Optimal Transport.
\end{IEEEbiography}

\begin{IEEEbiography}[{\includegraphics[width=1in,height=1.25in,clip,keepaspectratio, trim=0pt 23pt 0pt 0pt]{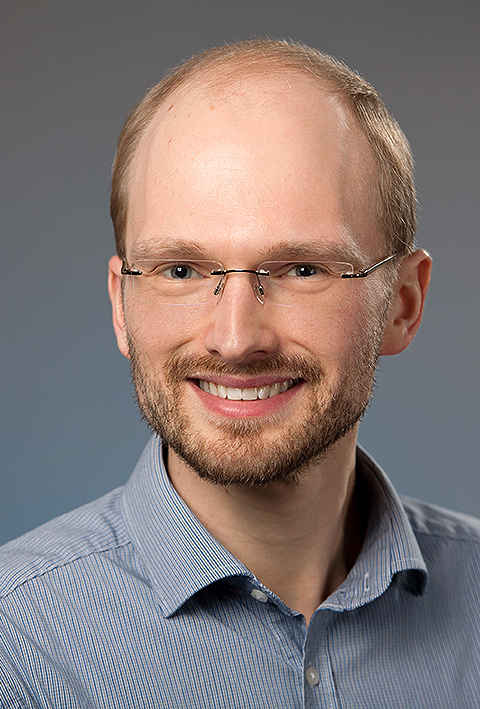}}]{Robert Beinert}
received his PhD in Mathematics from the University of
Göttingen (Germany).  After a Postdoc at the University of Graz
(Austria), he is currently a research fellow at the TU Berlin
(Germany).  His research interests include Inverse Problems, 
Optimization, Harmonic
Analysis, and Convex Analysis with applications in Signal and Image Processing
like Phase Retrieval.
\end{IEEEbiography}
 
\vspace{11pt}
  
\begin{IEEEbiography}[{\includegraphics[width=1in,height=1.25in,clip,keepaspectratio, trim=110pt 0pt 110pt 0pt]{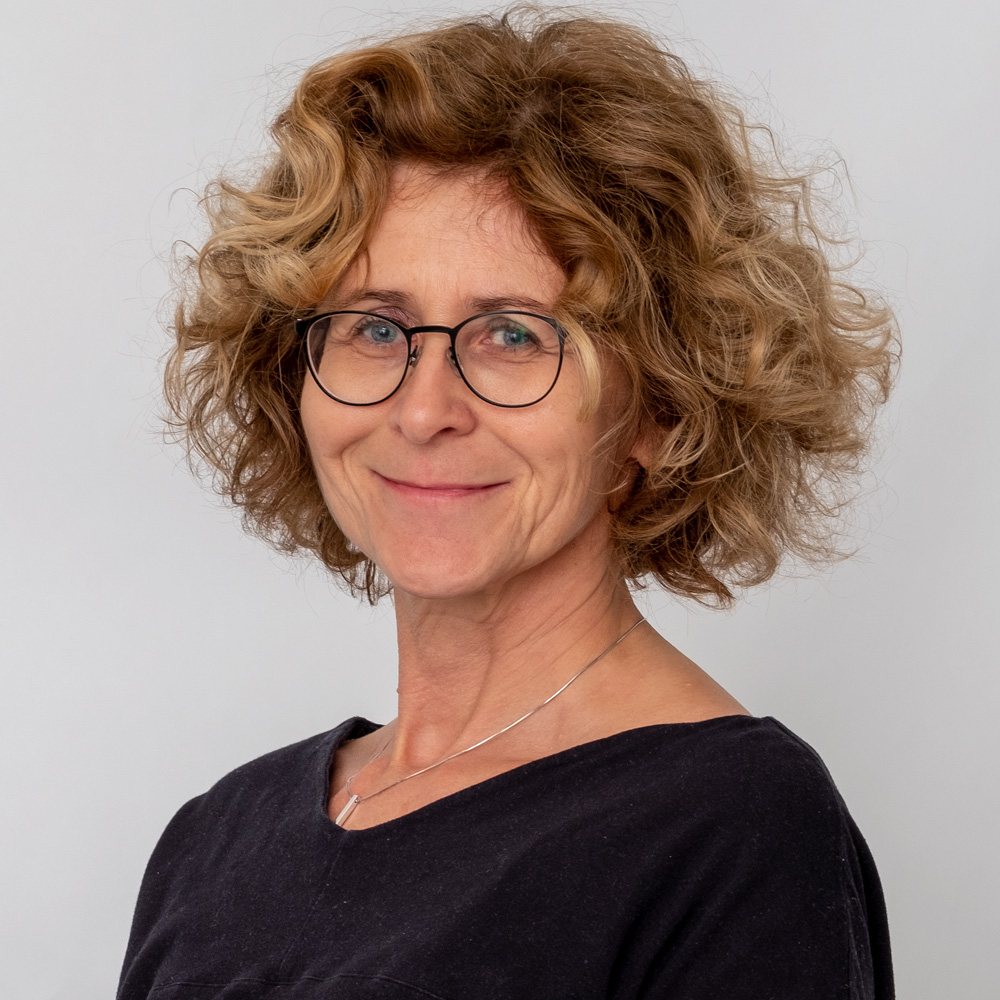}}]{Gabriele Steidl}
received her PhD and habilitation in mathematics from the University of Rostock.
After positions as associated professor for Mathematics at the TU Darmstadt and full professor at the University of Mannheim 
and TU Kaiserslautern, she is currently professor at the TU Berlin. 
She was a Postdoc, resp. visiting Professor at the Univ. of Debrecen, Zürich, ENS Cachan/Paris Univ. Paris Est, Sorbonne/IHP Paris 
and worked as consultant of the Fraunhofer ITWM Kaiserslautern. 
Gabriele Steidl is Editor-in-Chief of the SIAM Journal of Imaging Sciences and SIAM Fellow.
Her research interests include Harmonic Analysis, Optimization, Inverse Problem and Machine Learning
with applications in Image Processing.
\end{IEEEbiography}

\vfill
\end{document}